\documentclass[12pt]{amsart}
\usepackage{amsmath,amsfonts,amssymb,amscd,amsthm,amsbsy,epsf,graphicx, xcolor}
\usepackage{amsrefs}

\numberwithin{equation}{section}

\newtheorem{theorem}{Theorem}[section]
\newtheorem{lemma}[theorem]{Lemma}%[section]
\newtheorem{corollary}[theorem]{Corollary}
\newtheorem{proposition}[theorem]{Proposition}
\theoremstyle{definition}

\newtheorem{remark}[theorem]{Remark}
\def\bye{\end{document}} \def\by{\end{proof}\end{document}}

\def\R{{\mathbb R}}
\def\N{{\mathbb N}}

\def\mid{\,:\,}

\def\USC{\operatorname{USC}}
\def\LSC{\operatorname{LSC}}

\DeclareMathOperator{\AC}{AC}
\DeclareMathOperator{\cI}{\mathcal{I}}
\DeclareMathOperator{\cS}{\mathcal{S}}
\def\1{\mathbf{1}}

\def\Lip{\operatorname{Lip}}

\def\beq{\begin{equation}}
\def\eeq{\end{equation}}
\def\bald{\begin{aligned}}
\def\eald{\end{aligned}}

\def\cD{\mathcal{D}}
\def\ep{\varepsilon}

\def\N{\mathbb{N}}

\def\Lip{\operatorname{Lip}}
\def\T{\mathbb{T}}\def\R{\mathbb{R}}
\def\IN{\text{ in }}
\def\AND{\text{ and }}
\def\FORALL{\text{ for all }}
\def\FOR{\text{ for }}
\def\ON{\text{ on }}
\def\IN{\text{ in }}
\def\IF{\text{ if }}
\def\OW{\text{ otherwise}}
\def\bproof{\begin{proof}}
\def\eproof{\end{proof}}
\def\ga{\alpha}
\def\gth{\theta}
\def\ol{\overline}
  \newlength\ovlheight
       
\def\tim{\times}
\def\erf{\eqref}
\def\gl{\lambda}
\def\fr{\frac}
\def\stm{\setminus}
\def\bcases{\begin{cases}}
\def\ecases{\end{cases}}
\def\IF{\text{ if }}

\def\bye{\end{document}}

\def\gd{\delta}
\def\pl{\partial}
\def\du#1{\langle #1\rangle}
\def\gL{\Lambda}

\def\ctM{T^*\! M}
\def\gs{\sigma}
\def\cC{\mathcal{C}} \def\gam{\gamma}
\def\wtil{\widetilde}
\def\gz{\zeta}
\def\go{\omega} 
\def\BJ{\cS_{\mathrm{BJ}}}
\def\tH{\widetilde H} \def\hH{\widehat H}

\title[Hamilton-Jacobi equations]{
Hamilton-Jacobi equations with their Hamiltonians depending Lipschitz continuously on the unknown}

\author[Hitoshi Ishii, Kaizhi Wang, Lin Wang, and Jun Yan]{Hitoshi Ishii${}^*$ \and Kaizhi Wang \and Lin Wang \and Jun Yan}

\address[H. Ishii]{ Institute for Mathematics and Computer Science, Tsuda University, Tsuda 2-1-1, Kodaira, Tokyo 169-8050, Japan. }
\email{hitoshi.ishii@waseda.jp}
\address[K. Wang]{School of Mathematical Sciences, Shanghai Jiao Tong University, Shanghai 200240, China.}
\email{kzwang@sjtu.edu.cn}
\address[L. Wang]{School of Mathematics and Statistics, Beijing Institute of Technology, Beijing 100081, China.}
\email{lwang@bit.edu.cn}
\address[J. Yan]{School of Mathematical Sciences, Fudan University, Shanghai 200433, China.}\email{yanjun@fudan.edu.cn}


\thanks{H. Ishii was partially supported by the JSPS KAKENHI Grant Nos. JP16H03948, JP20K03688,
and JP20H01817; {K. Wang was partially supported by NSFC Grant Nos. 11771283, 11931016;
L. Wang was partially supported by NSFC Grant Nos.  11790273, 11631006;
J. Yan was partially  supported by NSFC Grant Nos.  11631006, 11790273.}}

\thanks{${}^*\,$Corresponding author}


\date{\today}

\begin{document}
\addtolength{\baselineskip}{3pt}

\maketitle
\begin{abstract} We study the Hamilton-Jacobi equations $H(x,Du,u)=0$ in $M$
and $\pl u/\pl t +H(x,D_xu,u)=0$ in $M\tim(0,\infty)$,
where the Hamiltonian $H=H(x,p,u)$ depends Lipschitz continuously on the variable $u$. In the framework of the semicontinuous viscosity solutions due to Barron-Jensen,
we establish the comparison principle, existence theorem, and representation formula
as value functions for extended real-valued, lower semicontinuous solutions for the Cauchy problem. We also establish some results on the long-time behavior of solutions for the Cauchy problem and classification of solutions for the stationary
problem.
 \end{abstract}

\tableofcontents

\section{Introduction}\label{sec:int}

We study the Hamilton-Jacobi equation
\beq\label{int.eq1}
u_t+H(x,D_xu,u)=0 \ \ \IN M\tim(0,T),
\eeq	
where $M$ is a connected, closed, smooth Riemannian manifold of dimension $n$,
$T$ is either a positive number or $+\infty$,
$H$ is the Hamiltonian on $\ctM\tim\R$,
$u$ is the unknown function on $M\tim[0,T)$, $u_t$ denotes the partial derivative $\pl u/\pl t=\pl_t u$, $D=D_x$ denotes the differential map, so that $(x, Du)=(x,Du(x))$ denotes an element of the cotangent space $T_{x}^* M$.
The Riemannian structure on $M$ induces a norm
$|\cdot|=|\cdot|_x$ on the tangent space $T_x M$. The canonical pairing between $T^*_x M$ and $T_x M$ is denoted by $\du{\cdot,\cdot}=\du{\cdot,\cdot}_x$, which defines naturally a
norm $|\cdot|=|\cdot|_x$ on $T^*_x M$.
%%%%%%%%%%%%%%%%%
\renewcommand{\theenumi}{H\arabic{enumi}}
\renewcommand{\labelenumi}{(\theenumi)}
%%%%%%%%%%%%%%%%%%%%%%

The following list collects our main {assumptions} on the Hamiltonian $H$.
\begin{enumerate}
\item%[(H1)]
\label{H1}The function $(x,p,u)\mapsto H(x,p,u)$ is continuous on $T^*M\tim\R$.
\item%[(H2)]
\label{H2}For any $R>0$ there exists $K>0$ such that
\[
H(x,p,u)>R \ \ \IF |u|\leq R \AND |p|\geq K.
\]

\item%[(H3)]
\label{H3}For any $(x,u)\in M\tim \R$, the function $p\mapsto H(x,p,u)$
is convex on $T_x^*M$.
\item %[(H4)]
\label{H4} The functions $u\mapsto H(x,p,u)$ are equi-Lipschitz continuous
on $\R$, with $(x,p)\in\ctM$.
\end{enumerate}
When \erf{H4} is assumed, the symbol $\gL$ denotes a positive Lipschitz bound:
\[
|H(x,p,u)-H(x,p,v)|\leq \gL|u-v| \ \ \FORALL (x,p)\in T^*M,\, u,v\in\R,
\]
and it is fixed throughout this paper.  Remark also that under the assumption
\erf{H4}, condition \erf{H2} is equivalent to say that for any $R>0$, there exists $K>0$ such that \ $H(x,p,0)>R$ \ if \ $|p|\geq K$.

In recent work K. Wang-L. Wang-J. Yan \cite{WWY2017}, the authors have studied the Cauchy problem
for \erf{int.eq1} with the initial condition of the form, with given $y\in M$ and $c\in\R$,
\beq\label{int.eq1.1}
u(x,0)=\bcases
c & \IF x=y \\
\infty& \OW.
\ecases
\eeq
To assure the existence of a solution belonging in $C(M\tim(0,T),\R)$, they assume a coercivity assumption stronger than (H2) above. 
The purpose of this paper is to adopt the notion of semicontinuous
viscosity solutions due to Barron-Jensen \cite{BJ1990} and to give an existence
 and uniqueness result for the Cauchy problem for \erf{int.eq1} with general
lower semicontinuous data. Consider the Hamilton-Jacobi equation
\beq\label{int.eq1.2}
u_t+|D_xu|-u=0  \ \ \IN M\tim(0,\infty),
\eeq
with the initial condition \erf{int.eq1.1}. This equation has the property of a finite speed propagation, by which any solution has necessarily discontinuities for positive times as far as the initial data is discontinuous. Indeed, the main concern of the paper
\cite{BJ1990} is to take into this kind of singular behavior of the solutions into
the notion of solution, and the solution $u$ of \erf{int.eq1.2} and \erf{int.eq1.1} is
given by the formula (see the proof of Proposition \ref{two-bounds} below
for a related discussion)
\[
u(x,t)=
\bcases
c\,e^t  \ \ & \ \IF d(x,y)\leq t,\\
+\infty & \ \OW. 
\ecases
\]

New features of our {{adaptation}} of the semicontinuous viscosity solutions
 to \erf{int.eq1} are the use of structure conditions on \erf{int.eq1}, that is,
the hypotheses \erf{H1}--\erf{H4} different {from} those employed in \cite{B1993, BJ1990, BJ1991,
F1993},  and the extension of the notion of solution which applies to extended
real-value functions.

We recall the definition of lower semicontinuous viscosity sub and supersolutions of
\erf{int.eq1} following \cite{I2001, BJ1990, F1993}: let $u: M\tim(0,T) \to \R\cup\{\infty\}$
be a lower semicontinuous function, that is, $u\in\LSC(M\tim(0,T),\R\cup\{\infty\})$.
The function $u$ is called a lower semicontinuous viscosity subsolution (resp., supersolution) of
\erf{int.eq1} if, whenever $\phi\in C^1(M\tim(0,T),\R)$ and
$\min(u-\phi)=(u-\phi)(x,t)$ for some $(x,t)\in M\tim(0,T)$, we have
$\phi_t(x,t)+H(x, D\phi(x,t),u(x,t))\leq 0$ (resp., $\geq 0$). If $u$ is both
lower semicontinuous sub and supersolutions of \erf{int.eq1},
then we call it a lower semicontinuous solution of \erf{int.eq1}.
Henceforth, for simplicity of notation, we write simply ``BJ" for ``lower semicontinuous
viscosity".  For instance, we say a BJ subsolution instead of a lower semicontinuous
viscosity subsolution.

We remark that a function $u\in \LSC(M\tim(0,T),\R\cup\{\infty\})$ is a BJ subsolution
of $u_t+H(x,D_xu,u)=0$ in $M\tim(0,T)$ if and only if it is a viscosity
supersolution, in the sense of Crandall-Lions (\cite{CL1983, CEL1984}), of
$-u_t-H(x,D_xu,u)=0$ in $M\tim(0,T)$.

Let $u\in \LSC(M\tim (0,T),\R\cup\{\infty\})$ and $(x,t)\in M\tim(0,T)$ be such that
$u(x,t)<\infty$. Let $D^-u(x,t)$ denote the subdifferential of $u$ at $(x,t)$ defined
as the set of all $(p,q)\in {T^*_x M}\tim\R$ such that in local coordinates, as $|y|+|s|\to 0$,
\[
u(x+y,t+s)\geq u(x,t)+\du{p,y}+qs+o(|y|+|s|).
\]

It is straightforward to generalize
the notion of BJ solution to a general Hamilton-Jacobi equation $F(x,Du,u)=0$ in $U$.
We henceforth write $\BJ(F)=\BJ(F,U)$ (resp., $\BJ^+(F)=\BJ^+(F,U)$ and
$\BJ^-(F)=\BJ^-(F,U)$) for the set of all BJ solutions (resp., BJ supersolutions and
BJ subsolutions) $u\in\LSC(U,\R\cup\{\infty\})$ of $F(x,Du,u)=0$ in $U$. For instance, $\BJ(\pl_t+H)=\BJ(\pl_t+H,M\tim(0,T))$ denotes the set of all BJ solutions of \erf{int.eq1}.
Similarly, we write  $\cS(F)=\cS(F,U)$ (resp., $\cS^+(F)=\cS^+(F,U)$ and
$\cS^-(F)=\cS^-(F,U)$) for the set of all viscosity solutions
$u\in\LSC(U,\R)$
(resp., viscosity  supersolutions $u\in\LSC(U,\R\cup\{\infty\})$ and
viscosity subsolutions $u\in\USC(U,\R\cup\{-\infty\})$ of $F(x,Du,u)=0$ in $U$
in the Crandall-Lions sense.

The remark after the introduction of BJ solutions above, can be stated as $\BJ^-(F)=\cS^+(-F)$.

Another simple remark here is that for any continuous function $u$,
we have $u\in\cS^-(F)$ if and only if $-u\in \BJ^-(F^\ominus)$, where
$F^\ominus$ is given by $F^\ominus(x,p,u)=F(x,-p,-u)$.

After this introduction, we establish a comparison principle for BJ sub and supersolutions of \erf{int.eq1} in Section 2, an existence result for the Cauchy
problem for \erf{int.eq1} is proved in Section 3, and a representation formula
for BJ solutions of \erf{int.eq1}, based on the idea of value functions of
optimal control associated with \erf{int.eq1} is presented in Section 4.
In Section 5, we introduce the notion of fundamental solutions
to \erf{int.eq1} as well as their basic properties. In Section 6,
we investigate the long-time behavior of solutions of \erf{int.eq1},
which is applied in Section 7 to classification of solutions to the corresponding
stationary problem together with several suggestive examples. The appendix presents a proof of a classical existence
result of Lipschitz continuous solutions to the Cauchy problem for  \erf{int.eq1}.

\section{Comparison principle} \label{comparison}

In this work, it is of major importance to establish the comparison principle for BJ solutions of \erf{int.eq1}.

\begin{theorem} \label{cwch} Assume \erf{H1}--\erf{H4}.  Let $v, w\in\LSC(M\tim[0,T),\R\cup\{\infty\})$ be, respectively, BJ sub and supersolutions of
\beq %\tag{\ref{ncwch}.1}
\label{comparison.1}
u_t+H(x,D_xu,u)=0 \ \ \IN \ M\tim(0,T).
\eeq
Assume that
\beq%\tag{\ref{ncwch}.2}
\label{comparison.2}
w(x,0) \geq \liminf_{t\to 0+}v(x,t)  \ \ \FORALL \ x\in M.
\eeq
Then, $v\leq w$ on $M\tim[0,T)$.
\end{theorem}

A stronger inequality than \erf{comparison.2} implies a stronger conclusion in the above theorem,
as stated in the following corollary.

\begin{corollary}\label{cor-cwch} Under the hypotheses of Theorem \ref{cwch},  but with \erf{comparison.2} replaced by the condition that
for some constant $C>0$,
\[
w(x,0)\geq C+\liminf_{t\to 0+}v(x,t) \ \ \FORALL x\in M,
\]
we have \ $v(x,t)+C e^{-\gL t}\leq w(x,t)\,$ for all $\,(x,t)\in M\tim[0,T)$.
Similarly, if \erf{comparison.2} is replaced by the condition that
for some constant $C>0$,
\[
w(x,0)+C\geq \liminf_{t\to 0+}v(x,t) \ \ \FORALL x\in M,
\]
then we have  \ $v(x,t)\leq w(x,t)+C e^{\gL t}\,$ for all $\,(x,t)\in M\tim[0,T)$.
\end{corollary}

\bproof Set $z(x,t)=v(x,t)+Ce^{-\gL t}$ and compute in a slightly informal fashion that
\[
z_t+H(x,D_xz,z)\leq v_t-\gL C e^{-\gL t}+H(x, D_xv,v)+\gL Ce^{-\gL t}\leq 0,
\]
to find that $z$ is a BJ subsolution of \erf{comparison.1}. It is then easily seen by Theorem \ref{cwch} that
$z\leq w$ in $M\tim(0,T)$. Instead, if we set $z(x,t)=w(x,t)+Ce^{\gL t}$, then
we find that $z$ is a BJ supersolution of \erf{comparison.1} and conclude by Theorem \ref{cwch} that $v\leq w$ in $M\tim[0,T)$.
\eproof

The next two lemmas constitute the primary part of the proof of Theorem \ref{cwch}.

\begin{lemma} \label{l1.cwch} %\label{l1.ncwch}
 In addition to the hypotheses of Theorem %\ref{ncwch}
\ref{cwch}, assume that
$v$ is Lipschitz continuous on $M\tim[0, T)$. Then $v\leq w$ on $M\tim[0,T)$.
\end{lemma}

It is a classical observation in the literature that the Lipschitz property of a viscosity subsolution or supersolution simplifies the formulation and proof of the comparison theorem. 

We remark that, under the hypotheses of Lemma \ref{l1.cwch}, condition \erf{comparison.2} is equivalent to the inequality
$v(x,0)\leq w(x,0)$ for all $x\in M$.

\bproof Thanks to the Lipschitz regularity of $v$, the function $v$ is a viscosity subsolution
of \erf{comparison.1} in the Crandall-Lions sense, which is a classical observation
due to \cite{BJ1990} (see also \cite[Theorem 2.3]{I2001}).
By the standard change of the unknown functions
(i.e., by considering the new unknowns $v(x,t)e^{-(\gL+1)t}$ and $w(x,t)e^{-(\gL+1)t}$),
we may assume that $u\mapsto H(x,p,u)-u$ is nondecreasing.

To the contrary to the conclusion, we suppose  that $\sup_{M\tim[0,T)}(v-w)>0$, and we will obtain a contradiction.
We can choose $S\in(0,T]$ so that $\sup_{M\tim[0, S)}(v-w)>0$.
Let $\ep>0$ and consider the function
\[
u(x,t)=v(x,t)-\fr{\ep}{S+\ep^2-t}  \ \ \ON M\tim[0,S].
\]
Note that, since $-w$ is bounded from above on $M\tim[0,S]$, if $\ep>0$ is sufficiently small, then
\[
\max_{x\in M}(u(x,S)-w(x,S))=\max_{x\in M}(v(x,S)-w(x,S))-\fr 1\ep<0.
\]
By assumption, we have
\[
\max_{x\in M}(u(x,0)-w(x,0))\leq -\fr{1}{S+\ep^2}<0.
\]
Choosing a point $(x_0,t_0)\in M\tim[0, S)$ such that
\[
(v-w)(x_0,t_0)>0,
\]
we observe that for sufficiently small $\ep>0$.
\[
(u-w)(x_0,t_0)=(v-w)(x_0,t_0)-\fr{\ep}{S+\ep^2-t_0}>0.
\]
Hence, fixing $\ep>0$ small enough, we find that
\[
\max_{M\tim[0,S]}(u-w)>0 \ \ \AND \ \ \max_{\pl (M\tim [0,S])}(u-w)<0.
\]
Note also that
\[\bald
u_t+H(x,D_xu,u)&=v_t-\fr{\ep}{(S+\ep^2-t)^2}+H\Big(x,D_xv,v-\fr{\ep}{S+\ep^2-t}\Big)
\\&\leq v_t+H(x,D_xv,v)\leq 0 \ \ \IN M \tim(0,S).
\eald\]

Fix a maximum point $(\hat x,\hat t)\in M\tim(0,S)$ of the function $u-w$ on
$M\tim[0,S]$. We choose a chart $(U,\phi)$ such that $\hat x\in U$.
We identify $U$ with $\phi(U)$, so that $\hat x\in U\subset \R^n$.
Consider the function
\[
\Phi_\ga: (x,t,y,s)\mapsto u(x,t)-w(y,s)-\ga |x-y|^2-\ga (t-s)^2-|x-\hat x|^2-(t-\hat t)^2
\]
on $U\tim[0,S]\tim U\tim [0,S]$. We fix a compact neighborhood $B\subset U \tim[0,S]$
of $(\hat x,\hat t)$.  Let $(x_\ga,t_\ga,y_\ga,s_\ga)\in B\tim B$ be a maximum point
of $\Phi_\ga$ on the set $B\tim B$.  It is a standard observation that as $\ga\to \infty$,
\[
(x_\ga,t_\ga,y_\ga,s_\ga) \to (\hat x,\hat t,\hat x,\hat t) \ \ \AND \ \ w(y_\ga,s_\ga) \to
w(\hat x,\hat t)
\]
Assuming $\ga$ large enough, we may assume that
$x_\ga,y_\ga$ are in the interior of $B$. By the viscosity properties of $u$ and $w$, we find that
\[\bald
&2\ga(t_\ga-s_\ga)+2(t_\ga-\hat t)+H(x_\ga,2\ga(x_\ga-y_\ga)+2(x_\ga-\hat x), u(x_\ga,t_\ga))\leq 0,
\\&2\ga(t_\ga-s_\ga)+H(y_\ga, 2\ga(x_\ga-y_\ga),w(y_\ga,s_\ga))\geq 0
\eald
\]
The Lipschitz continuity of $u$ implies that the collections
$\{\ga(x_\ga-y_\ga)\}$ and $\{\ga(t_\ga-s_\ga)\}$ are bounded in $\R^n$ and $\R$,
respectively.
We may choose a sequence $\{\ga_j\} \subset (0,\infty)$ so that, as $j\to \infty$,
\[
2\ga_j(x_{\ga_j}-y_{\ga_j}) \to \hat p \ \ \AND \ \ 2\ga_j(t_{\ga_j}-s_{\ga_j}) \to \hat q.
\]
Sending $j\to \infty$, we get
\[\bald
&\hat q+H(\hat x,\hat p, u(\hat x,\hat t))\leq 0,
\\&\hat q+H(\hat x, \hat p,w(\hat x,\hat t))\geq 0,
\eald
\]
and, subtracting one from the other and recalling that $u\mapsto H(x,p,u)-u$ is nondecreasing,
\[
0\geq H(\hat x,\hat p, u(\hat x,\hat t))-H(\hat x,\hat p, w(\hat x,\hat t))
\geq (u-w)(\hat x,\hat t)>0.
\]
This is a contradiction, which proves that $v\leq w$ on $M\tim[0,T)$.
\eproof

\begin{lemma} \label{l2.cwch} Assume \erf{H1}--\erf{H4}.  Let $v\in\LSC(M\tim[0,T),\R\cup\{\infty\})$ be a BJ subsolution of \ $
u_t+H(x,D_xu,u)=0 \ \IN M \tim(0,T)$,
and assume that $v$ is bounded from below on $M\tim[0,T)$.
Let $C_0, C_1>0$ be constants such that
\[
H(x,p,0)\geq -C_0 \ \ \FORALL (x,p)\in\ctM \ \ \AND \ \ v\geq -C_1 \ \ \ON M\tim[0,T).
\]
Fix $y\in M$. If \ $\liminf_{t\to 0+} v(y,t)<\infty$, then
\beq \label{comparison.3}
v(y,t)\leq e^{\gL t}\,\liminf_{t\to 0+}v(y,t) +(C_0\gL^{-1}+2C_1)(e^{\gL t}-1) \ \ \FORALL t\in (0,T),
\eeq
Furthermore, for any $s\in(0,T)$, if \ $v(y,s)<\infty$, then
\beq \label{comparison.4}
v(y,t)\leq v(y,s)e^{\gL (t-s)}+(C_0\gL^{-1}+2C_1)(e^{\gL (t-s)}-1) \ \ \FORALL t\in (s,T).
\eeq
\end{lemma}

We remark that if \ $u\in\LSC(M\tim[0,T),\R\cup\{\infty\})$, then
$u$ is bounded from below on $M\tim[0,S]$ for any $0<S<T$.

\bproof  We first show that \erf{comparison.3} is valid.
We set\ $C_2=C_0+2\gL C_1$\ and compute in an informal way that
\[
0\geq v_t+H(x,D_xv,v)\geq v_t+H(x,p,0)-\gL|v|
\geq v_t-C_0-\gL(v+2C_1),
\]
and deduce that $v$ is a subsolution of
\beq\label{l2.visco}
v_t-\gL v-C_2 = 0  \ \ \IN M\tim(0,T).
\eeq

Set $c=\liminf_{t\to 0+} v(y,t)$. Suppose, to the contrary to \erf{comparison.3}, that there is $\tau\in(0,T)$ such that
\[
v(y,\tau)>ce^{\gL \tau}+C_2\gL^{-1}(e^{\gL \tau}-1).
\]
Note that if we define the function $v_0$ on $[0,T)$ by
\[
v_0(t)=\bcases
c &\IF t=0, \\
v(y,t) &\OW
\ecases
\]
then $v_0\in \LSC([0,T),\R\cup\{\infty\})$.
Choose $\ep>\gd>0$ so small that
\[
v_0(\tau)>\ep+v_0(0)e^{\gL\tau} +(C_2+\ep)\gL^{-1}(e^{\gL\tau}-1)
\ \ \ \AND \ \ \
\gL\gd <\ep.
\]
For any $a\in\R$, set
\[
\psi_a(t)=ae^{\gL t}+(C_2+\ep)\gL^{-1}(e^{\gL t}-1),
\]
and note that
\[
\psi_a'(t)-\gL\psi_a(t)-C_2=\ep \ \ \FORALL t\geq 0.
\]
Fix a function $\chi \in C^1([0,\infty),\R)$ such that
\[
\chi(0)=0, \quad
\chi'(t)\geq 0 \ \ \FORALL t\geq 0, \ \ \AND \ \ \chi(t)=\gd \ \  \FORALL t\geq \tau/2,
\]
and for any $k\in\N$, set
\[
\chi_k(t)=\chi(kt) \ \ \AND \ \ f_{a,k}(t)=\psi_a(t)+\chi_k(t) \ \ \FOR t\geq 0.
\]
Observe that for all $t\geq 0$,
\beq\label{comparison.5}
f_{a,k}'(t)-\gL f_{a,k}(t)-C_2= \ep +\chi_k'(t)-\gL \chi_k(t)
\geq \ep -\gL\gd>0.
\eeq

Consider the function $f_k:=\psi_{a,k}$ with $a:=v_0(0)$.
Since \ $v_0(0)=c=\liminf_{t\to 0+} v_0(t)$,
we may choose $s\in(0,\tau)$ such that
\ $v_0(s)<v_0(0)e^{\gL s}+\gd$.
If $k\in\N$ is so large that $ks>\tau/2$, then \ $\chi_k(s)=\gd$,
and hence,
\[
f_k(s)=v_0(0)e^{\gL s}+(C_2+\ep)(e^{\gL s}-1)+\gd>v_0(s).
\]
Fix such $k\in\N$ and note that
\[
\min_{t\in[0,\tau]}(v_0-f_k)(t)\leq (v_0-f_k)(s)<0,
\]
and, by the choice of $\tau,\,\ep,\,\gd$,
\[
(v_0-f_k)(0)=0,\qquad (v_0-f_k)(\tau)>0.
\]
It is clear that the function
\[
\Phi: a\mapsto \min_{t\in[0,\tau]}(v_0-f_{a,k})(t)
\]
is monotone decreasing, continuous,  and satisfies
\[
\lim_{a\to\infty}\Phi(a)=-\infty \ \ \AND \ \ \lim_{a\to -\infty}\Phi(a)=\infty.
\]
These observations show that there is a unique $a<v_0(0)$ such that
\[
\min_{t\in[0,\tau]}(v_0-f_{a,k})(t)=0,
\]
while
\[
\min_{t=0, \tau}(v_0-f_{a,k})(t)>\min_{t=0,\tau}(v_0-f_k)(t)=0.
\]
Thus, there is $t_0\in(0,\tau)$ such that $t\mapsto v_0(t)-f_{a,k}(t)$ takes
a minimum $0$ at $t=t_0$.  Since $f_{a,k}(t)=v_0(t)=v(y,t)$ at $t=t_0$,  by
 \erf{comparison.5} we have
\beq \label{comparison.6}
f_{a,k}'(t_0)-\gL v(y,t_0)-C_2=f_{a,k}'(t_0)-\gL f_{a,k}(t_0)-C_2>0.
\eeq

Let $a,k$ be as above and set $f=f_{a,k}$ for simplicity of notation.
Fix $\gs>0$ so that $\gs<t_0<\tau$, and consider the function
\[
\Psi_\ga(x,t):= v(x,t)-f(t)+(t-t_0)^2 + \ga d(x,y)^2
\]
on the set $M\tim [\gs,\tau]$, where $\ga>0$ is a constant to be sent to $\infty$.
Let $(x_\ga,t_\ga)$ be a minimum point of $\Psi_\ga$. It is easily seen that
$(x_\ga,t_\ga) \to (y,t_0)$ and  $v(x_\ga, t_\ga) \to v(y,t_0)$ as $\ga\to \infty$. Thus, if $\ga$ is large enough,
$(x_\ga,t_\ga)$ is an interior point of $M\tim[\gs,\tau]$ and hence, by \erf{l2.visco},
\[
f'(t_\ga)-2(t_\ga-t_0)-\gL v(x_\ga,t_\ga)-C_0\leq 0.
\]
Sending $\ga\to\infty$ yields
\[
f_{a,k}'(t_0)-\gL v(y,t_0)-C_0\leq 0,
\]
which contradicts \erf{comparison.6} and completes the proof of \erf{comparison.3}.

Now, we prove \erf{comparison.4}.  Fix any $s\in(0,T)$ such that $v(y,s)\in\R$,
and set
\[
\eta(t)=v(y,s)e^{\gL(t-s)}+C_2\gL^{-1}(e^{\gL(t-s)}-1) \ \ \FOR t\in[0,T).
\]
If \ $
\liminf_{t\to 0+} v(y,t)\leq \eta(0)$,
then we find by \erf{comparison.3} that for all $t\in[0,T)$,
\[
v(y,t)\leq \eta(0)e^{\gL t}+C_2\gL^{-1}(e^{\gL t}-1)
=v(y,s)e^{\gL(t-s)}+C_2\gL^{-1}(e^{\gL(t-s)}-1),
\]
which shows that \erf{comparison.4} holds.

Next, consider the case \ $\liminf_{t\to 0+}v(y,t)>\eta(0)$.
As before, define $v_0\in\LSC([0,T),\R\cup\{\infty\})$ by
\[
v_0(t)=
\bcases
\displaystyle
\liminf_{t\to 0+}v(y,t) &\IF t=0, \\[2pt] v(y,t) & \IF t>0.
\ecases
\]
Notice that the value $v_0(0)$ can be $+\infty$.

To the contrary to \erf{comparison.4}, we suppose that for some $\tau\in(s,T)$,
\[
v_0(\tau)> \eta(\tau),
\]
and will get a contradiction. Note that
\[
v_0(0)>\eta(0), \quad v_0(\tau)>\eta(\tau), \ \ \AND \ \  v_0(s)=\eta(s).
\]
We choose $\ep>0$ so that
\[
v_0(\tau)>\eta(\tau)+\ep C_2\gL^{-1}(e^{\gL(\tau-s)}-1),
\]
and, for $a\in\R$, define the smooth function $g_a$ on $[0,T)$ by
\[
g_a(t)
=a e^{\gL(t-s)}+(C_2+\ep)\gL^{-1}(e^{\gL(t-s)}-1).
\]
Observe that if $a=v_0(s)$, then
$g_a(t)=\eta(t)+\ep C_2\gL^{-1}(e^{\gL(t-s)}-1)$ for all $t\in [0,T)$ and
\beq\label{comparison.7}
v_0(0)>g_a(0), \quad v_0(\tau)>g_a(\tau), \ \ \AND \ \  v_0(s)=g_a(s),
\eeq
which implies
\[
\min_{[0,\tau]}(v_0-g_a)\leq 0.
\]
Note also that, as $a\to-\infty$, $\min_{[0,\tau]}(v_0-g_a)$ goes
monotonically and continuously to $+\infty$. Consequently, there is $a\leq v_0(s)$
such that
\[
\min_{[0,\tau]}(v_0-g_a)=0,
\]
while, by \erf{comparison.7}, $(v_0-g_a)(t)>0$ for $t=0,\tau$.  Hence, there
is $t_0\in(0,\tau)$ such that $t\mapsto v(y,t)-g_a(t)$ takes a local minimum at
$t=t_0$, with value $0$.
As in the proof of \erf{comparison.3} above, \erf{l2.visco} and
the above observation yield
\[
g_{a}'(t_0)-\gL g_{a}(t_0)-C_2=g_{a}'(t_0)-\gL v(y,t_0)-C_2\leq 0,
\]
while the function $g_{a}$ satisfies
\[
g_{a}'(t)-\gL g_{a}(t)-C_2=\ep \ \ \FORALL t>0.
\]
This contradiction completes the proof. \eproof

The following lemma is a standard observation about BJ subsolutions.

\begin{lemma} \label{BJ-min} Let $u,v\in\BJ^-(F, U)$, where
$(F,U)$ is either $(H,M)$ or $(\pl_t+H,M\tim(0,T))$.
Then, $\min\{u,v\}\in\BJ^-(F,U)$. If $v,w\in\BJ(F,U)$ instead, then
$\min\{v,w\}\in\BJ(F,U)$.
\end{lemma}

\bproof Let $\psi\in C^1(U,\R)$ and $z\in U$ be a minimum point of the function
$\min\{u,v\}-\psi$. Observe that if $\min\{u,v\}=u$ at $z$, then $u-\psi$ takes a minimum at $z$, and otherwise $v-\psi$ takes a minimum at $z$. If
$u,v\in \BJ^-(F,U)$ (resp., $u,v\in \BJ(F,U)$), we find that $F(z,D\psi(z),\min\{u,v\}(z))\leq 0$ (resp., $F(z,D\psi(z),\min\{u,v\}(z))=0$), which completes the proof.
\eproof

The following proof of Theorem \ref{cwch} has a strong similarity to that of 
\cite[Theorem 4.1]{I2008} (see also \cite[Theorem 1]{I1997}).

\bproof[Proof of Theorem \ref{cwch}]  We need only to prove
that $v\leq w$ on $M\tim(0,S)$ for all $0<S<T$. Thus, we may assume
in this proof that $T<\infty$ and $v$ is bounded from below on $M\tim[0,T)$.

We first observe that we may assume that $v$ is a bounded function on $M\tim[0,T)$.
Choosing $C>0$ so large that
\[
H(x,0,0)\leq C\gL e^{-2\gL T},
\]
and setting $z(x,t)=C e^{-2\gL t}$ for $(x,t)\in M\tim[0,T)$, we compute
\[
z_t+H(x,D_x z,z)\leq -2\gL z+H(x,0,0)+\gL z \leq 0 \ \ \ON M\tim[0,T),
\]
to find that $z$ is a subsolution of \erf{comparison.1}. Thanks to
Lemma \ref{BJ-min}, the function $\min\{v,z\}$ is a BJ subsolution of \erf{comparison.1}, which is bounded on $M\tim[0,T)$. It is obvious that \erf{comparison.2} is valid with $\min\{v,z\}$ in place of $v$.
To conclude that $v\leq w$ on $M\tim[0,T)$, we only need to show that  $\min\{v,z\}\leq w$ on $M\tim[0,T)$ for all $C$ large enough.
Thus, we may assume by replacing $v$ by $\min\{v,z\}$, with $C$ sufficiently large, if necessary that $v$ is bounded on $M\tim[0,T)$.

{Next}, we regularize $v$ by the inf-convolution in the variable $t$. We define the function $v_0$ on $M\tim[0,T)$ by
\[
v_0(x,t)=\bcases\displaystyle
\liminf_{t\to 0+}v(x,t) & \text{ for } \  t=0,\\[2pt]
v(x,t) & \text{ for }\  t>0,
\ecases
\]
let $\ep>0$, and set
\[
v_\ep(x,t)=\inf_{s\in[0,T)} \Big(v_0(x,s)+\fr{e^{-\gL t}}{2\ep}(t-s)^2\Big) \ \ \FOR (x,t)\in M\tim [0,T).
\]

Fix a constant $C_0>0$ such that $H(x,p,0)\geq -C_0$ for all $(x,p)\in T^*M$ and
$|v(x,t)|\leq C_0$ for all $(x,t)\in M\tim[0,T)$.  Using \erf{comparison.3} in Lemma \ref{l2.cwch}, we obtain
\beq \label{comparison.10}
v_0(x,t)\leq v_0(x,0) +C_0(1+\gL)(e^{\gL t}-1)
\leq v_0(x,0)+C_1t,
\eeq
where $C_1:=C_0 e^{\gL T}\gL(1+\gL)$.

We set
\[
T_1(\ep)=2\ep C_1 e^{\gL T} \quad \AND\quad T_2(\ep)=T-2\sqrt{\ep(C_0+1)e^{\gL T}},
\]
assume that $\ep>0$ is sufficiently small so that $T_1(\ep)<T_2(\ep)$,
and prove that
\beq\label{comparison.11}
v_\ep(x,t)=\min_{s\in (0,T)}\Big(v_0(x,s)+\fr{e^{-\gL t}}{2\ep}(t-s)^2\Big) \ \ \FORALL (x,t)\in M\tim[T_1(\ep), T_2(\ep)],
\eeq
that is, the minimum above is attained at some $s\in[0,T)$.

We first prove that
\beq\label{comparison.12}
v_\ep(x,t)=\min_{s\in [0,T)}\Big(v_0(x,s)+\fr{e^{-\gL t}}{2\ep}(t-s)^2\Big) \ \ \FORALL (x,t)\in M\tim[0, T_2(\ep)],
\eeq
To see this, we fix any $(x,t)\in M\tim[0,T_2(\ep)]$.
If $s\in(0,T)$ is chosen so that $v_\ep(x,t)+1> v_0(x,s)+e^{-\gL t}(t-s)^2/(2\ep)$, then
\[
C_0+1\geq v_\ep(x,t)+1 \geq -C_0+\fr{e^{-\gL t}}{2\ep}(t-s)^2,
\]
and, consequently,
\[
|t-s|< 2\sqrt{\ep(C_0+\tfrac12)e^{\gL T}}\ \ \AND \ \ s< T_2(\ep)+2\sqrt{\ep(C_0+\tfrac 12)e^{\gL T}}<T.
\]
Applying this estimate to any minimizing  sequence
$\{s_k\}_{k\in\N}\subset(0,T)$ of the minimization problem in the definition of $v_\ep(x,t)$, we find that
\erf{comparison.12} holds. Also, applying the estimates above to the minimizer $s$
in \erf{comparison.12}, we find that for all $(x,t)\in M\tim[0,T_2(\ep)]$,
\beq\label{comparison.13}
v_\ep(x,t)\geq \inf\left\{v_0(x,s)\mid 0\leq s<T,\ |s-t|<2\sqrt{\ep(C_0+1)e^{\gL T}}\right\}.
\eeq

To complete the proof of \erf{comparison.11}, we fix $(x,t)\in M\tim[T_1(\ep),T_2(\ep)]$
and note by \erf{comparison.10} and \erf{comparison.12} that
\beq\label{comparison.14}
v_\ep(x,t)\leq v_0(x,t)\leq v_0(x,0)+C_1 t.
\eeq
Hence, if $s=0$ were a minimizer of \erf{comparison.12}, then, by the choice of $T_1(\ep)$, we would have
\[\bald
v_0(x,0)+C_1 t&\geq
v_\ep(x,t)=v_0(x,0)+\fr{e^{-\gL t}}{2\ep}t^2
\\&> v_0(x,0)+\fr{T_1(\ep)e^{-\gL T}}{2\ep}t=v_0(x,0)+C_1 t,
\eald
\]
which is a contradiction. This together with \erf{comparison.12} assures that
\erf{comparison.11} is valid.

By \erf{comparison.14} and \erf{comparison.2}, we have for all $x\in M$,
\beq \label{comparison.15}
v_\ep(x,T_1(\ep))\leq v_0(x,0)+C_1T_1(\ep)\leq w(x,0)+C_1T_1(\ep).
\eeq
By the definition of $v_\ep$, the family of the functions $t\mapsto v_\ep(x,t)$, with $x\in M$, is equi-Lipschitz continuous on $[0,T)$.

Now, we claim that $v_\ep$ satisfies, in the BJ sense,
\[
v_{\ep,t}+H(x,D_xv_{\ep},v_\ep) \leq 0 \ \ \IN M\tim(T_1(\ep),T_2(\ep)).
\]
To check this, let $\phi\in C^1(M\tim(T_1(\ep),T_2(\ep)),\R)$ and assume that
$v_\ep-\phi$ takes a minimum at $(\hat x,\hat t)\in M\tim(T_1(\ep),T_2(\ep))$
and for some $\hat s\in(0,T)$,
\[
v_\ep(\hat x,\hat t)=v_0(\hat x,\hat s)+\fr{e^{-\gL \hat t}}{2\ep}(\hat t-\hat s)^2.
\]
Then, the function
\[
(x,t,s) \mapsto v_0(x,s)+\fr{e^{-\gL t}}{2\ep}(t-s)^2-\phi(x,t)
\]
takes a (local) minimum at $(\hat x,\hat t,\hat s)$. This implies that
\[\bald
&\Big(D_x\phi(\hat x,\hat t), \fr{e^{-\gL\hat t}}{\ep}(\hat t-\hat s)\Big)\in D^-v_0(\hat x,\hat s),
\\& \phi_t(\hat x,\hat t)=\fr{e^{-\gL \hat t}}{\ep}(\hat t-\hat s)-\fr{\gL\,e^{-\gL \hat t}}{2\ep}(\hat t-\hat s)^2.
\eald\]
Since $v_0$ is a BJ subsolution of $u_t+H(x,D_xu,u)=0$ in $M\tim(0,T)$, we have
\[\bald
0&\geq \fr{e^{-\gL \hat t}}{\ep}(\hat t-\hat s)+H(\hat x,D_x\phi(\hat x,\hat t),v_0(\hat x,\hat s))
\\&=\phi_t(\hat x,\hat t)+\gL\fr{e^{-\gL \hat t}}{2\ep}(\hat t-\hat s)^2
+H\Big(\hat x,D_x\phi(\hat x,\hat t),v_\ep(\hat x,\hat t)-\fr{e^{-\gL \hat t}}{2\ep}(\hat t-\hat s)^2\Big)
\\&\geq \phi_t(\hat x,\hat t)+\gL\fr{e^{-\gL \hat t}}{2\ep}(\hat t-\hat s)^2
+H(\hat x,D_x\phi(\hat x,\hat t),v_\ep(\hat x,\hat t))-\fr{\gL\,e^{-\gL \hat t}}{2\ep}(\hat t-\hat s)^2
\\&=\phi_t(\hat x,\hat t)+H(\hat x,D_x\phi(\hat x,\hat t),v_\ep(\hat x,\hat t)),
\eald
\]
which assures that $v_\ep$ is a BJ subsolution of $u_t+H(x,D_xu,u)=0$ in $M\tim(T_1(\ep),T_2(\ep))$.

Recalling that the functions $t\mapsto v_\ep(x,t)$, with $x\in M$, are equi-Lipschitz continuous on $[T_1(\ep),T_2(\ep)]$, we choose a constant $C(\ep)>0$
as a Lipschitz bound of the functions above, we find that $v_\ep$ is a BJ subsolution of
\[
-C(\ep)+H(x,D_xu, u)\leq 0 \ \ \IN M\tim(T_1(\ep),T_2(\ep)).
\]
Thus, since $v_\ep$ is a bounded function, we deduce that $v_\ep\in \Lip(M\tim [T_1(\ep),T_2(\ep)],\R)$.

Noting that the functions $(x,t)\mapsto v_\ep(x,t+T_1(\ep))-C_1T_1(\ep) e^{\gL t}$
is a subsolution of $u_t+H(x,D_xu,u)=0$ in $M\tim(0,T_2(\ep)-T_1(\ep))$ and recalling
\erf{comparison.15}, we invoke Lemma \ref{l1.cwch}, to obtain
\[
w(x,t)\geq v_\ep(x,t+T_1(\ep)) -C_1 T_1(\ep) e^{\gL  t} \ \ \FORALL (x,t)\in M\tim[0,T_2(\ep)-T_1(\ep)).
\]
Combine this with  \erf{comparison.13}, to get for all $(x,t)\in M\tim[0,T_2(\ep)-T_1(\ep))$,
\[\bald
w(x,t)&+C_1T_1(\ep) e^{\gL t}\geq
v_\ep(x,t+T_1(\ep))
\\& \geq \inf \left\{v_0(x,s)\mid s\in[0,T),\, |s-t|\leq T_1(\ep)+2\sqrt{\ep(C_0+1)e^{\gL T}}\right\},
\eald\]
which yields, in view of the lower semicontinuity of $t\mapsto v_0(x,t)$, that
for all $(x,t)\in M\tim[0,T)$,
\[
v(x,t)\leq v_0(x,t)\leq \liminf_{[0,T)\ni s\to t}v_0(x,s)\leq \liminf_{\ep\to 0+}v_\ep(x,t+T_1(\ep))\leq w(x,t).
\]
This completes the proof.
\eproof

\section{Existence of BJ solutions} \label{cp}

We consider the Cauchy problem
\beq\label{cp.1}\left\{\bald
&u_t+H(x,D_xu,u)=0 \ \ \IN M\tim[0,\infty),\\
&u(x,0)=\phi(x) \ \ \FOR x\in M,
\eald\right.
\eeq
where
\beq\label{cp.2}
\phi\in\LSC(M,\R\cup\{\infty\}).
\eeq
In view of Theorem \ref{cwch}, we understand the initial condition in \erf{cp.1} as
\[
\phi(x)=\liminf_{t\to 0+}u(x,t) \ \ \FORALL x\in M.
\]
We call a function $u\in\LSC(M\tim[0,\infty),\R\cup\{\infty\})$ a BJ solution
of \erf{cp.1} if $u$ is a BJ solution of $u_t+H(x,D_xu,u)=0$ in $M\tim(0,\infty)$
and
\beq \label{cp.3}
\phi(x)=u(x,0)=\liminf_{t\to \infty}u(x,t) \ \ \FORALL x\in M.
\eeq

\begin{theorem} \label{thm1.cp}Assume \erf{H1}--\erf{H4} and \erf{cp.2}. Then there exists
a BJ solution of \erf{cp.1}.
\end{theorem}

This theorem and Theorem \ref{cwch} assure that for each $\phi\in \LSC(M,\R\cup\{\infty\})$, there exists a unique BJ solution $u$ of \erf{cp.1}. In what follows
we write $S_t=S_t^H$ for the map $\phi\mapsto u(\cdot,t)$ for all $t\geq 0$.
It follows from Theorem \ref{cwch} and Lemma \ref{l2.cwch} that
for every $s,t\geq 0$, $S_t\circ S_s=S_{t+s}$.

\bproof We fix a constant $C_0>0$ so that
\[
H(x,p,0)\geq -C_0 \ \ \FORALL (x,p) \in T^*M \ \ \AND \ \ \phi\geq -C_0 \ \ \ON M.
\]
Select a sequence $\{\phi_k\}_{k\in\N}${$\subset C^1(M,\mathbb{R})$}
 such that for all $x\in M$,
\[
-C_0\leq \phi_k(x)\leq \phi_{k+1}(x) \ \ \AND \ \ \lim_{k\to\infty}\phi_k(x)=\phi(x).
\]
It is well-known (see Theorem \ref{a.1} in the appendix) that for each $k\in\N$ there is a Crandall-Lions
viscosity solution $u_k\in \Lip(M\tim[0,\infty),\R)$ of \erf{cp.1}, with $\phi_k$ in place of $\phi$, which is aslo
a BJ solultion of
\beq\label{cp.1.1}
u_t+H(x,D_xu,u)=0 \ \ \IN M\tim (0,\infty).
\eeq
Note that if we set $v(x,t)=-C_0 e^{\gL t}$
for $(x,t)\in M\tim[0,\infty)$, then $v$ is a subsolution of \erf{cp.1.1}.

By the classical comparison theorem (or Theorem \ref{cwch} above), we find that for all $k\in\N$,
\[
-C_0e^{\gL t}\leq u_k(x,t)\leq u_{k+1}(x,t) \ \ \FORALL (x,t)\in M\tim[0,\infty).
\]
We define a function $u_\infty$ on $M\tim[0,\infty)$ by
\[
u_\infty(x,t)=\lim_{k\to\infty} u_k(x,t).
\]
It is obvious that $u_\infty\in\LSC(M\tim[0,\infty),\R\cup\{\infty\})$.

We claim that $u_\infty$ is a BJ solution of \erf{cp.1}.
By the standard stability property of viscosity solutions, we find that $u_\infty$
is a BJ solution of $u_t+H(x,D_xu,u)=0$ in $M\tim(0,\infty)$. Indeed, as we see below,
$u_\infty$ is the so-called lower relaxed limit of $\{u_k\}$. It is clear that
for any $(x,t)\in M\tim(0,T)$,
\beq\label{cp.1.1+}
u_\infty(x,t)\geq \sup_{j\in\N}\inf\{u_k(y,s)\mid d(y,x)+|s-t|<j^{-1},\, k\geq j\}.
\eeq
On the other hand, for any $a<u_\infty(x,t)$, by the continuity of the $u_k$,
we may choose $l,m\in\N$ such that
\[
a<u_l(y,s) \ \ \IF \ d(y,x)+|s-t|<m^{-1},
\]
which implies that for $j:=\max\{l,m\}$,
\[
a<u_k(y,s) \ \ \IF \ k\geq j, \ d(y,x)+|s-t|<j^{-1}.
\]
This combined with \erf{cp.1.1+} shows that for all $(x,t)\in M\tim(0,T)$,
\beq\label{cp.1.1++}
u_\infty(x,t)=\sup_{j\in\N}\inf\{u_k(y,s)\mid d(y,x)+|s-t|<j^{-1},\, k\geq j\},
\eeq
where the formula on the right hand side is the lower relaxed limit of $\{u_k\}$.
The identity above can be stated as
\beq\label{cp.1.1+++}
-u_\infty(x,t)=\inf_{j\in\N}\sup\{-u_k(y,s)\mid d(y,x)+|s-t|<j^{-1},\, k\geq j\},
\eeq
where the formula on the right hand side is the upper relaxed limit of $\{-u_k\}$.
Noting that $-u_k\in\cS(-(\pl_t+H)^\ominus, M\tim(0,T))
\subset \cS^-((\pl_t+H)^\ominus, M\tim(0,T))$, we find
by \cite[Lemma 6.1]{CIL1992}, for instance,
that $-u_\infty\in \cS^-((\pl_t+H)^\ominus,M\tim(0,T))$, which implies that
$u_\infty\in \BJ^-(\pl_t+H,M\tim(0,T))$. It is straightforward from \erf{cp.1.1++} and
\cite[Remark 6.2]{CIL1992} to find that $u_\infty\in \cS^+(\pl_t+H, M\tim(0,T))$. These inclusions together prove that $u_\infty$ is a BJ solution of  \erf{cp.1.1}.

It remains to check \erf{cp.3}. Fix any $x\in M$. Note first that
\[
u_\infty(x,0)=\lim_{k\to\infty}u_k(x,0)=\lim_{k\to\infty}\phi_k(x)=\phi(x).
\]
For any $k\in\N$, we have
\[
\liminf_{t\to 0+} u_\infty(x,t)\geq \lim_{t\to 0+}u_k(x,t)=\phi_k(x),
\]
and hence,
\[
\phi(x)\leq \liminf_{t\to 0+}u_\infty(x,t).
\]
Now, \erf{comparison.3} assures that for all $0<T<\infty$ and
$(x,t)\in M\tim[0,T)$,
\[
u_k(x,t)\leq u_k(x,0)e^{\gL t}+C_0(\gL^{-1}+2)e^{\gL T}(e^{\gL t}-1) \ \
\FORALL k\in\N,
\]
which implies
\[
u_\infty(x,t)\leq \phi(x)e^{\gL t}+C_0(\gL^{-1}+2)e^{\gL T}(e^{\gL t}-1),
\]
and moreover,
\[
\liminf_{t\to 0+}u_\infty(x,t)\leq \phi(x)=u_\infty(x,0) \ \  \FORALL
x\in M.
\]
Thus, we find that \erf{cp.3} holds and $u_\infty$ is a BJ solution of \erf{cp.1}.
\eproof

\begin{corollary} \label{cor1.cp}
Under the hypotheses of Theorem \ref{thm1.cp}, there exists a unique $BJ$ solution
of \erf{cp.1}.
\end{corollary}

\bproof The existence and uniqueness assertions follow from Theorems \ref{thm1.cp}
and \ref{cwch}, respectively.
\eproof

\section{Value function representation}

We give a value function representation or the Hopf-Lax-Oleink formula for the solution of \erf{cp.1}.

Let $L$ denote the Lagrangian associated with $H$, that is,
\[
L(x,\xi,u)=\sup_{p\in T_x^*M} (\du{p,\xi}-H(x,p,u)) \ \ \FOR (x,\xi,u)\in TM\tim\R.
\]
Note that $L$ is lower semicontinuous on $TM\tim\R$, and that $(x,\xi,u)$ is coercive
and, furthermore, has a superlinear growth in $\xi$. To see the superlinear growth, let
$A>0$ be any constant and observe that
\beq\label{repr.super}
L(x,\xi,u)\geq \max_{p\in T^*_x M,\,|p|=A}(\du{p,\xi}-C(A,x,u))\geq A|\xi|-C(A,x,u),
\eeq
where {$C(A,x,u):=\max_{p\in T_x^*M,\, |p|=A}H(x,p,u)$}.

Let $u\in \LSC(M\tim[0,\infty),\R\cup\{+\infty\})$.
Note that for each $T>0$, $u$ is bounded from below on $M\tim[0,T]$
by a constant. For $(x,t)\in M\tim(0,\infty)$, let $\cC(x,t,u)$
denote the set of all $\gam\in\AC([0,t],M)$ such that
$\gam(t)=x$ and
\beq\label{curveC}
\int_0^t(|L(\gam(s),\dot\gam(s),0)|+|u(\gam(s),s)|)ds<\infty.
\eeq

\begin{theorem} \label{thm:repr}Assume \erf{H1}--\erf{H4}. Let $\phi\in\LSC(M,\R\cup\{\infty\})$ and
let $u$ be the BJ solution of \erf{cp.1}.  Fix $(x,t)\in M\tim(0,\infty)$
so that $u(x,t)<\infty$. Then
\beq\label{u-repr}
u(x,t)=\min_{\gam\in\cC(x,t,u)} \int_0^t L(\gam(s),\dot\gam(s),u(\gam(s),s))ds+\phi(\gam(0)),
\eeq
and the minimum above is attained.
\end{theorem}

It is convenient to convert our Hamilton-Jacobi equation to the one whose
Lagrangian $\wtil L(x,\xi,t,u)$ is increasing in $u$. Let $u$ be as in Theorem \ref{thm:repr}
and for $\gl\in\R$, set
\[
v(x,t)=e^{\gl t}u(x,t),
\]
and calculate  in a slightly informal way that
\[
v_t=e^{\gl t}(u_t+\gl u)=-e^{\gl t}H(x,D_xu,u)+\gl v=-e^{\gl t}H(x,e^{-\gl t} D_xv, e^{-\gl t}v)
+\gl v,
\]
to find that in the viscosity sense,
\[
v_t+\wtil H(x,D_xv,t,v)=0 \ \ \IN M\tim(0,\infty),
\]
where $\wtil H$ is given by
\[
\wtil H(x,p,t,u)=e^{\gl t}H(x,e^{-\gl t}p,e^{-\gl t}u)-\gl u \ \ \FOR (x,p,t,u)\in T^*M
\tim \R^2.
\]
The Lagrangian $\wtil L(x,\xi,t,u)$ corresponding to $\wtil H$ is defined as
\[\bald
\wtil L(x,\xi,t,u):&=\sup_{p\in T^*_xM}\du{p,\xi}-\wtil H(x,p,t,u)
\\&=e^{\gl t}L(x,\xi,e^{-\gl t}u)+\gl u \ \ \FOR (x,\xi,t,u)\in TM\tim\R^2.
\eald\]

Note that if $\gl\geq\gL$, then $u\mapsto \wtil L(x,\xi,t,u)$ is
nondecreasing on $\R$.  Henceforth, we fix
\[
\gl=\gL+1.
\]
As in the proof of Theorem \ref{thm1.cp}, we select a sequence
$\{\phi_k\}\subset C^1(M,\R)$ such that for all $x\in M$,
\[
\phi_k(x)\leq \phi_{k+1}(x) \ \ \AND \ \
\lim_{k\to\infty}\phi_k(x)=\phi(x).
\]
Let $u_k$ be the solution of \erf{cp.1}, with $\phi_k$ in place of $\phi$.  Note that
$u_k\in \Lip(M\tim[0,\infty),\R)$, $u_k\leq u_{k+1}$ on $M$ for all $k\in\N$,
and
\[
u(x,t)=\lim_{k\to \infty} u_k(x,t) \ \ \FORALL (x,t)\in M\tim[0,\infty).
\]

Set $v_k(x,t)=e^{\gl t}u_k(x,t)$ for $(x,t)\in M\tim[0,\infty)$.
We recall that for every $k\in\N$,
\beq \label{v_k-repr}
v_k(x,t)=\min_{\gam \in\AC([0,t],M), \, \gam(t)=x} \int_0^t
\wtil L(\gam(s),\dot\gam(s),s,v_k(\gam(s),s))ds+\phi_k(\gam(0)),
\eeq
and the minimum above is achieved at some $\gam\in\AC([0,t],M)$
satisfying $\gam(t)=x$. Note that, in the above formula,
\beq\bald\label{L-Ltil}
\int_0^t
&\wtil L(\gam(s),\dot\gam(s),s,v_k(\gam(s),s))ds
\\&=\int_0^t e^{\gl s}[L(\gam(s),\dot\gam(s),u_k(\gam(s),s))+\gl u_k(\gam(s),s)]ds.
\eald\eeq

\begin{lemma} \label{equiv}Assume \erf{H1}--\erf{H4}. Let $(x,t)\in M\tim(0,\infty)$ and $\gam\in\AC([0,t],M)$ be such that
$\gam(t)=x$. Let $u\in\LSC(M\tim [0,t],\R\cup\{+\infty\})$.
Then, $\gam\in \cC(x,t,u)$ if and only if
\beq \label{finite}
\int_0^t e^{\gl s}[L(\gam(s),\dot\gam(s),u(\gam(s),s))+\gl u(\gam(s),s)]ds<\infty.
\eeq
\end{lemma}

Notice that if $u\geq -C_0$ for some $C_0>0$, then
\begin{gather*}
L(y,\xi,u)\geq L(y,\xi,-C_0)-\gl(u+C_0),\\
\intertext{and }
e^{\gl s}[L(y,\xi,u)+\gl u]
\geq e^{\gl s}[L(x,\xi,-C_0)+\gl(-C_0)] \geq -C
\end{gather*}
for all $(y,\xi)\in TM$ and some constant $C>0$. Hence, the condition
\erf{finite} makes sense.

\bproof We fix a constant $C_0>0$ so that for all $(y,s)\in M\tim[0,t]$ and
$\xi\in T_y M$,
\[
u(y,s)\geq -C_0 \ \ \AND \ \ L(y,\xi,0)\geq -C_0,
\]
which, in particular, yield
\[
|u(y,s)|\leq 2C_0+u(y,s) \ \ \AND \ \ |L(y,\xi,0)|\leq 2C_0+L(y,\xi,0).
\]

Assume first that $\gam\in \cC(x,t,u)$.
Note that for any $(y,\xi,u)\in TM\tim\R$,
\[
L(y,\xi,u)+\gl u\leq L(y,\xi,0)+\gL |u|+\gl u
\leq |L(y,\xi,0)|+2\gl|u|,
\]
and hence,
\[\bald
\int_0^T e^{\gl s}&[L(\gam(s),\dot\gam(s),u(\gam(s),s))+\gl u(\gam(s),s)]ds
\\&\leq \int_0^T e^{\gl s} (|L(\gam(s),\dot\gam(s),0)|+2\gl |u(\gam(s),s)|)ds
\\&\leq e^{\gl T}\int_0^T (|L(\gam(s),\dot\gam(s),0)|+2\gl |u(\gam(s),s)|)ds.
\eald
\]
This shows that \erf{finite} holds.

Next, assume that \erf{finite} is satisfied.
Note that for any $(y,\xi,u)\in TM\tim\R$, if $u\geq -C_0$, then
\[\bald
|L(y,\xi,0)|+|u|&\leq 2C_0+L(y,\xi,0)+|u|\leq
2C_0+L(y,\xi,u)+\gL |u|+|u|
\\&\leq 2C_0(1+\gl)+L(y,\xi,u)+\gl u.
\eald\]
Hence, we have
\[\bald
\int_0^t &(|L(\gam(s),\dot\gam(s),0)|+|u(\gam(s),s)|)ds
\leq \int_0^t e^{\gl s}[|L(\gam(s),\dot\gam(s),0)|+|u(\gam(s),s)|]ds
\\&\leq \int_0^t e^{\gl s}[L(\gam(s),\dot\gam(s),u(\gam(s),s))+\gl u(\gam(s),s)]ds
+2C_0(1+\gl) \gl^{-1}(e^{\gl t}-1),
\eald
\]
which shows that $\gam\in\cC(x,t,u)$. The proof is complete.
\eproof

\bproof[Proof of Theorem \ref{thm:repr}] By formula \erf{v_k-repr},
we have
\[
v(x,t)\geq v_k(x,t)=\min_{\gam(0)=x}\int_0^t\wtil L(\gam(s),\dot\gam(s),s,v_k(\gam(s),s))ds+\phi_k(\gam(0)).
\]
We select a minimizer $\gam_k\in \AC([0,t],M)$, with $\gam_k(t)=x$,
for each $k\in\N$ in the above, to obtain
\[\bald
v(x,t)&\geq\int_0^t\wtil L(\gam_k(s),\dot\gam_k(s),s,v_k(\gam_k(s),s)))ds+\phi_k(\gam_k(0))
\\&\geq \int_0^t\wtil L(\gam_k(s),\dot\gam_k(s),s,v_j(\gam_k(s),s))ds+\phi_j(\gam_k(0))
\ \IF k\geq j,
\eald\]
where we have used the fact that $u\mapsto \wtil L(y,\xi,s,u)$ is nondecreasing.
Since $\xi \mapsto \wtil L(y,\xi,s,u)$ has a superlinear growth (see \erf{repr.super}), we may
select a subsequence of $\{\gam_k\}$, which will be denoted again by the same symbol, such that  for some $\gam\in \AC([0,t],M)$, as $k\to\infty$,
\[
\gam_k \to \gam \ \IN C([0,t],M)
\ \ \AND \ \ \dot\gam_k \to \dot\gam
\ \hbox{ weakly in } L^1([0,t],TM).
\]
It follows that for any $j\in\N$,
\[
v(x,t)\geq \int_0^t\wtil L(\gam(s),\dot\gam(s),s,v_j(\gam(s),s))ds+\phi_j(\gam(0)).
\]
Furthermore, by the monotone convergence theorem,
\beq\label{min-for-v}
v(x,t)\geq \int_0^t\wtil L(\gam(s),\dot\gam(s),s,v(\gam(s),s))ds+\phi(\gam(0)).
\eeq
As noted in \erf{L-Ltil}, we have
\[
\int_0^t\wtil L(\gam(s),\dot\gam(s),s,v(\gam(s),s))ds
=\int_0^t e^{\gl s}[L(\gam(s),\dot\gam(s),u(\gam(s),s))+\gl u(\gam(s),s)]ds,
\]
and hence, \erf{min-for-v} assures together with Lemma \ref{equiv}
that $\gam\in \cC(x,t,u)$.

On the other hand, from \erf{v_k-repr}, we have for any $\eta\in\cC(x,t,u)$,
\[\bald
v_k(x,t)&\leq\int_0^t\wtil L(\eta(s),\dot\eta(s),s,v_k(\eta(s),s))ds+\phi_k(\eta(0))
\\&\leq\int_0^t\wtil L(\eta(s),\dot\eta(s),s,v(\eta(s),s))ds+\phi(\eta(0)),
\eald
\]
and moreover,
\[
v(x,t)\leq\int_0^t\wtil L(\eta(s),\dot\eta(s),s,v(\eta(s),s))ds+\phi(\eta(0)).
\]
Thus, we have
\beq\label{v-repr}
v(x,t)=\min_{\eta\in\cC(x,t,u)}\int_0^t \wtil L(\eta(s),\dot\eta(s),s,v(\eta(s),s))ds+\phi(\eta(0)).
\eeq

We now deduce from \erf{v-repr} that \erf{u-repr} is valid.

Fix any $\eta\in \cC(x,t,u)$.
We may assume that $\phi(\eta(0))<\infty$. Indeed, otherwise, it is obvious that
\[
u(x,t)\leq \int_0^t L(\eta(s),\dot\eta(s),u(\eta(s),s))ds+\phi(\eta(0)).
\]
Note that $\eta\in \cC(\eta(\tau),\tau,u)$ for all $\tau\in(0,t]$ and that
$u(\eta(s),s)<\infty$ for a.e. $s\in[0,t]$.
By \erf{v-repr}, which is valid for general $(x,t)$, we have for a.e. $\tau\in(0,t]$,
\[\bald
e^{\gl\tau}u(\eta(\tau),\tau)
&=v(\eta(\tau),\tau)\leq \int_0^\tau \wtil L(\eta(s),\dot\eta(s),s,v(\eta(s),s))ds+\phi(\eta(0))
\\&=\int_0^\tau e^{\gl s}[L(\eta(s),\dot\eta(s),u(\eta(s),s))+\gl u(\eta(s),s)]ds+\phi(\eta(0)).
\eald
\]
Setting
\[
f(\tau)=e^{-\gl\tau} \Big(\int_0^\tau e^{\gl s}[L(\eta(s),\dot\eta(s),u(\eta(s),s))+\gl u(\eta(s),s)]ds+\phi(\eta(0))\Big),
\]
and using the above, we compute that for a.e. $\tau\in(0,t)$,
\[
f'(\tau)=-\gl f(\tau)+
L(\eta(\tau),\dot\eta(\tau),u(\eta(\tau),\tau))+\gl u(\eta(\tau),\tau)
\leq L(\eta(\tau),\dot\eta(\tau),u(\eta(\tau),\tau)),
\]
and moreover,
\[\bald
u(x,t)&\leq f(t)\leq \int_0^t L(\eta(s),\dot\eta(s),u(\eta(s),s))ds+f(0)
\\&= \int_0^t L(\eta(s),\dot\eta(s),u(\eta(s),s))ds+\phi(\eta(0)).
\eald\]
Hence, we find that
\[
u(x,t)\leq \inf_{\eta\in\cC(x,t,u)}
\int_0^t L(\eta(s),\dot\eta(s),u(\eta(s),s))ds+\phi(\eta(0)).
\]

Now, let $\gam\in\cC(x,t,u)$ be a minimizer for the right hand side of \erf{v-repr}.
We claim that for a.e. $\tau\in(0,t)$,
\beq \label{v=int}
v(\gam(\tau),\tau)=\int_0^\tau\wtil L(\gam(s),\dot\gam(s),s,v(\gam(s),s))ds+
\phi(\gam(0)).
\eeq
Once this is proved, setting
\[
f(\tau)=e^{-\gl\tau}\left(\int_0^\tau e^{\gl s}[L(\gam(s),\dot\gam(s),u(\gam(s),s))+\gl u(\gam(s),s)]ds+
\phi(\gam(0))\right),
\]
we argue as above, to find that for a.e. $\tau\in(0,t)$,
\[
f'(\tau)=L(\gam(\tau),\dot\gam(\tau),u(\gam(\tau),\tau))
\]
and
\[
u(x,t)=f(t) =\int_0^t L(\gam(s),\dot\gam(s),u(\gam(s),s))ds+\phi(\gam(0)),
\]
which shows that \erf{u-repr} is valid and there exists a minimizer of the
minimization in \erf{u-repr}.

It remains to prove that \erf{v=int} holds.
As before, from \erf{v-repr} we find that for a.e. $\tau\in(0,t)$,
\[
v(\gam(\tau),\tau)\leq \int_0^\tau\wtil L(\gam(s),\dot\gam(s),s,v(\gam(s),s))ds+
\phi(\gam(0)).
\]
We only need to show that the above inequalities are in fact equalities.
To see this, we suppose to the contrary that for some $\tau\in(0,t)$,
\beq \label{strict-ineq}
v(\gam(\tau),\tau)< \int_0^\tau\wtil L(\gam(s),\dot\gam(s),s,v(\gam(s),s))ds+
\phi(\gam(0)).
\eeq
By \erf{v-repr}, there is a curve $\eta\in\cC(\gam(\tau),\tau,u)$ such that
\[
v(\gam(\tau),\tau)=\int_0^\tau \wtil L(\eta(s),\dot\eta(s),s,v(\eta(s),s))ds+\phi(\eta(0)).
\]
We define a curve $\gz\in\AC([0,t],M)$ by setting
\[
\gz(s)=\bcases
\eta(s) & \FOR s\in[0,\tau], \\
\gam(s) & \FOR s\in(\tau,t],
\ecases
\]
and note that $\gz\in\cC(x,t,u)$. Observe by \erf{v-repr} and \erf{strict-ineq} that
\[\bald
v(x,t)&\leq
\int_0^t \wtil L(\gz(s),\dot\gz(s),s,v(\gz(s), s))ds+\phi(\gz(0))
\\&=\int_\tau^t \wtil L(\gam(s),\dot\gam(s),s,v(\gam(s), s))ds
+\int_0^\tau \wtil L(\eta(s),\dot\eta(s),s,v(\eta(s),s))ds+\phi(\eta(0))
\\&=\int_\tau^t \wtil L(\gam(s),\dot\gam(s),s,v(\gam(s), s))ds
+v(\gam(\tau),\tau)
\\&<\int_\tau^t \wtil L(\gam(s),\dot\gam(s),s,v(\gam(s), s))ds
+\int_0^\tau \wtil L(\gam(s),\dot\gam(s),s,v(\gam(s),s))ds+\phi(\gam(0))
\\&=v(x,t),
\eald
\]
which is a contradiction. This shows that \erf{v=int} holds.
\eproof

\section{Fundamental solultions}
Given $\phi\in\LSC(M,\R\cup\{\infty\})$,
if $u=u(x,t)$ is the BJ solution of \erf{cp.1}, then we write $u(x,t,\phi):=u(x,t)$ for notational clarity.
Let $(y,c)\in M\tim (\R\cup\{\infty\})$ and define $\phi\in\LSC(M,\R\cup\{\infty\})$ by
\beq \label{f1}
\phi_{y,c}(x)=\bcases
c&\IF x=y,\\[2pt]
\infty &\OW.
\ecases
\eeq
We set $h(x,t,y,c)=u(x,t,\phi_{y,c})$ for $(x,t)\in M\tim[0,T)$. Notice that $h(x,t,y,\infty)=\infty$.
We call this function $h(x,t,y,c)$ on $M\tim[0,T)$, with parameter
$(y,c)\in M\tim (\R\cup\{\infty\})$, a \emph{fundamental solution} to $u_t+H(x,D_xu,u)=0$ in $M\tim(0,T)$.

\begin{lemma} \label{f1.lemma} Assume \erf{H1}--\erf{H4}.
\begin{enumerate}
\item[(i)]For any $(x,t,y)\in M\tim[0,T)\tim M$, the function
$c\mapsto h(x,t,y,c)$ is {nondecreasing} on $\R$ and Lipschitz continuous on $\R$, with Lipschitz bound $e^{\gL t}$.
\item[(ii)] The function $h$ is lower semicontinuous on $M\tim[0,T)\tim M\tim {(\R\cup\{\infty\})}$.
\item[(iii)]For any $\phi\in\LSC(M,\R\cup\{\infty\})$, the function $(x,t,y)\mapsto h(x,t,y,\phi(y))$ is lower semicontinuous on $M\tim [0,T)\tim M$.
\item[(iv)] For any $\phi\in\LSC(M,\R\cup\{\infty\})$, the function $(x,t)\mapsto \inf_{y\in M}h(x,t,y,\phi(y))$ is lower semicontinuous on $M\tim[0,T)$.
\end{enumerate}
\end{lemma}

Before going into the proof of the above lemma, we recall that, by definition,
a neighborhood of $\infty\in (-\infty,\infty]=\R\cup\{\infty\}$ is a subset
of $(-\infty,\infty]$ containing a set of the form $(a,\infty]=(a,\infty)\cup\{\infty\}$,
with $a\in\R$.

\bproof We begin with assertion (i). Let $(x,t,y)\in M\tim[0,T)\tim M$ and $c_1,c_2\in\R$. Let
$\phi_{y,c_1}$ and $\phi_{y,c_2}$ be the functions defined by \erf{f1}, with
$c=c_1, c_2$, respectively. If $c_1<c_2$, then $\phi_{y,c_1}<\phi_{y,c_2}$ on $M$,
and Theorem \ref{cwch} yields $h(x,t,y,c_1)\geq h(x,t,y,c_2)$ for all
$(x,t,y)\in M\tim[0,T)\tim M$. This assures the desired monotonicity of $h(x,t,y,c)$ in $c$. Setting $v(x,t)=h(x,t,y,c_1)+|c_1-c_2|e^{\gL t}$, we note that $v$ is a BJ
supersolution of $u_t+H(x,D_xu,u)=0$ in $M\tim(0,T)$ and that $h(x,0,y,c_2)\leq v(x,0)$ for all $x\in M$, and conclude by Theorem \ref{cwch} that
$h(x,t,y,c_2)\leq h(x,t,y,c_1)+|c_1-c_2|e^{\gL t}$ for all $(x,t,y)\in M\tim[0,T)\tim M$.
This shows that $c\mapsto h(x,t,y,c)$ is Lipschitz continuous with Lipschitz bound $e^{\gL t}$ for any $(x,t,y)\in M\tim[0,T)\tim M$.

To check (ii), let $(x_0,t_0,y_0,c_0)\in M\tim [0,T)\tim M\tim (\R\cup\{\infty\})$
and $a\in\R$ be such that
$h(x_0,t_0,y_0,c_0)>a$.  By assertion (i), there is $b<c_0$ such that
$h(x_0,t_0,y_0,b)>a$.
In view of the proof of
Theorem \ref{thm1.cp}, choosing a sequence $\{\phi_k\}_{k\in\N}\subset
C^1(M,\R)$ such that
\[
\phi_k\leq \phi_{y_0, b}\ \ \ON M \ \ \AND \ \ \lim_{k\to\infty}\phi_k(x)=\phi_{y_0,b}(x)
\ \ \FORALL x\in M,
\]
we have
\[
\lim_{k\to\infty} u(x,t,\phi_k)=h(x,t,y_0,b) \ \ \FORALL (x,t)\in M\tim[0,T).
\]
Fix $k\in\N$ so that
\[
u(x_0,t_0,\phi_k)>a.
\]
By the continuity of $(x,t)\mapsto u(x,t,\phi_k)$ and $\phi_k$, we may select neighborhoods $V$ and $W$ of $(x_0,t_0)$ and $(y_0,c_0)$, respectively, so that
\[
u(x,t,\phi_k)>a \ \ \FORALL (x,t)\in V \quad\AND\quad
\phi_k\leq \phi_{y,c} \ \ \ON M \ \ \FORALL (y,c)\in W,
\]
which imply, together with Theorem \ref{cwch}, that
\[
h(x,t,y,c)\geq u(x,t,\phi_k)>a \ \ \FORALL (x,t,y,c)\in V\tim W.
\]
This assures that (ii) is valid.

To prove assertion (iii), we fix $\phi\in \LSC(M,\R\cup\{\infty\})$.
Let $(x_0,t_0,y_0)\in M\tim [0,T)\tim M$
and $a\in\R$ be such that $h(x_0,t_0,y_0,\phi(y_0)>a$.  By assertion (ii), we can choose a {neighborhood} $V$ of $(x_0,t_0,y_0,\phi(y_0))$ such that
\[
h(x,t,y,c)>a \ \ \FORALL (x,t,y,c)\in V.
\]
In view of the monotonicity of $c\mapsto h(x,t,y,c)$, we may assume that
$V=W\tim(b,\infty]$ for some neighborhood $W$ of $(x_0,t_0,y_0)$.
By the semicontinuity of $\phi$, we can choose a {neighborhood} $U$ of $y_0$
so that $\phi(y)>b$ for all $y\in U$. Then, we have
\[
h(x,t,y,\phi(y))>a \ \ \FORALL (x,t,y)\in W\cap (M\tim [0,T)\tim U),
\]
which shows the lower semicontinuity of $(x,t,y)\mapsto h(x,t,y,\phi(y))$, proving (iii).

Now, we prove (iv). Fix any $\phi\in\LSC(M,\R\cup\{\infty\})$.
We note that  for any $(x,t)\in M\tim[0,T)$,
the function $y\mapsto h(x,t,y,\phi(y))$ is lower semicontinuous on $M$
by assertion (iii), $M$ is compact, and therefore, it attains a minimum at some point
in $M$.
Let $(x_0,t_0)\in M\tim[0, T)$ and $a\in \R$ be such that
$\min_{y\in M}h(x_0,t_0,y,\phi(y))>a$. Noting that
\[
h(x_0,t_0,y,\phi(y))>a \ \ \FORALL y\in M,
\]
in view of assertion (iii), for each $y\in M$ we can choose
neighborhoods $U_y$ and $V_y$ of $(x_0,t_0)$ and $y$, respectively,
such that
\[
h(x,t,z,\phi(z))>a \ \ \FORALL (x,t,z)\in U_y\tim V_y.
\]
Since $M$ is compact, we may select a finite number of $y_i\in M
$, with $i=1,\ldots,J$, so that $M=\bigcup_{i=1}^J V_{y_i}$. Then, we set
$U=\bigcap_{i=1}^J U_{y_i}$, to find that
\[
h(x,t,z,\phi(z))>a \ \ \FORALL (x,t,z)\in U\tim M.
\]
This shows that for the {neighborhood} $U$ of $(x_0,t_0)$ and all $(x,t)\in U$,
\[\min_{y\in M}h(x,t,y,\phi(y))>a,\] which proves assertion (iv).
\eproof

\begin{theorem} \label{f2.thm}Assume \erf{H1}--\erf{H4}. Let $\phi\in\LSC(M,\R\cup\{\infty\})$
and let $u\in \LSC(M\tim[0,T),\R\cup\{\infty\})$ be the (unique) BJ solution of \erf{cp.1}.
Then,
\[
u(x,t)=\min_{y\in M}h(x,t,y,\phi(y)) \ \ \FORALL (x,t)\in M\tim[0,T).
\]
\end{theorem}

\bproof Set
\[
v(x,t)=\min_{y\in M}h(x,t,y,\phi(y)) \ \ \FOR (x,t)\in M\tim[0,T).
\]
Lemma \ref{f1.lemma}, (iv) {assures} that $v\in\LSC(M\tim[0,T),\R\cup\{\infty\})$.
It is a standard observation that $v$ is a BJ solution of $u_t+H(x,D_xu,u)=0$
in $M\tim(0,T)$. Indeed, let $\psi\in C^2(M\tim(0,T),\R)$ and assume that
$v-\psi$ takes a {minimum} at $(x_0,t_0)\in M\tim(0,T)$. Choose $y_0\in M$ so that
\[
v(x_0,t_0)=h(x_0,t_0,y_0,\phi(y_0)),
\]
and note that $(x,t)\mapsto h(x,t,y_0,\phi(y_0))-\psi(x,t)$ takes a minimum
at $(x_0,t_0)$, which yields, since $h$ is a fundamental solution,
\[\bald
0&=\psi_t(x_0,t_0)+H(x_0,D_x\psi(x_0,t_0), h(x_0,t_0,y_0,\phi(y_0)))
\\&=\psi_t(x_0,t_0)+H(x_0,D_x\psi(x_0,t_0), v(x_0,t_0)).
\eald
\]
Hence, $v$ is a BJ solution of $u_t+H(x,D_xu,u)=0$
in $M\tim(0,T)$.

For each fixed $y\in M$, we have $\phi(x)\leq \phi_{y,\phi(y)}(x)$ for all $x\in M$,
and moreover,  by Theorem \ref{cwch}, $u(x,t)\leq h(x,t,y,\phi(y))$ for all
$(x,t)\in M\tim[0,T)$. Since $y\in M$ is arbitrary, we find that
$u(x,t)\leq v(x,t)$ for all $(x,t)\in M\tim[0,T)$. For fixed $y\in M$, if $\phi(y)<\infty$,
then, by the definition of fundamental solutions (see also \erf{cp.3}),
\[
\phi(y)=\liminf_{t\to 0+}h(y,t, y,\phi(y)),
\]
which implies
\[
u(y,0)=\phi(y)\geq \liminf_{t\to 0+}v(y,t).
\]
Thus, we have
\[
u(x,0)\geq \liminf_{t\to 0+}v(x,t) \ \ \FORALL x\in M.
\]
Applying Theorem \ref{cwch}, we find that $u(x,t)\geq v(x,t)$ for all $(x,t)
\in M\tim[0,T)$. Thus, we have $u=v$ on $M\tim[0,T)$.
\eproof

In our generality, the fundamental solution $h(x,t,y,c)$ may take value $+\infty$
at some point $(x,t)\in M\tim(0,\infty)$. A simple example is as follows.
Consider the case where $M=\T^1$ and $H(x,p,u)=|p|$. In this case we have
\[
h(x,t,y,c)=\bcases
c &\IF d(x,y)\leq t,\\
+\infty & \OW.
\ecases
\]

Under hypotheses \erf{H1} and \erf{H2}, there are constants $C_0>0$ and $R>0$ such that
for all $(x,p)\in T^*M$,
\[
H(x,0,0)\leq C_0 \quad \AND \quad H(x,p,0)\geq C_0+1\ \ \IF \ |p|\geq R.
\]
Moreover, if (H3) holds, then for all $(x,p)\in T^*M$,
\[
H(x,p,0)\geq \tfrac 1 R |p|+C_0 \ \ \IF\ |p|\geq R.
\]
It is now obvious that if \erf{H1}--\erf{H4} hold, then there exist constants $\gd>0$
and $C_1>0$ such that
\beq\label{linear-growth}
H(x,p,u)\geq \gd|p|-C_1-\gL|u| \ \ \FORALL (x,p,u)\in T^*M\tim\R.
\eeq

\begin{proposition} \label{two-bounds} Assume \erf{H1}--\erf{H4}. Let
$y\in M$. There exist constants $\gd>0$ and $C>0$
 such that for all $(x,t)\in M\tim[0,\infty)$,
\begin{align}
h(x,t,y,0)&\geq -C(e^{\gL t}-1),  \label{l-bound}
\\ h(x,t,y,0)&\leq C(e^{\gL t}-1) \ \ \IF \ d(x,y)\leq \gd t, \label{u-bound}
\end{align}
\end{proposition}

\bproof  To check \erf{l-bound}, fix a constant $C_0>0$ so that
\[
H(x,0,0)\leq C_0 \ \ \FORALL x\in M,
\]
set
\[
v(x,t)=-\gL^{-1}C_0(e^{\gL t}-1) \ \ \FOR (x,t)\in M\tim[0,\infty),
\]
and note that $v$ is a subsolution of $u_t+H(x,D_xu,u)=0$
in $M\tim(0,\infty)$.
Since $v\in C(M\tim[0,\infty),\R)$ and $v(x,0)=0\leq h(x,0,y,0)$ for all $x\in M$,
Theorem \ref{cwch} yields that
\[
h(x,t,y,0)\geq v(x,t)=-\gL^{-1}C_0(e^{\gL t}-1) \ \ \FORALL (x,t)\in M\tim[0,\infty).
\]

To show \erf{u-bound}, fix constants $\gd>0$ and $C_1>0$ so that \erf{linear-growth} holds and define $w\in \LSC(M\tim[0,\infty),\R\cup\{\infty\})$ by
\[
w(x,t)=
\bcases
\gL^{-1}C_1(e^{\gL t}-1) & \IF \ d(x,y)\leq \gd t,\\
+\infty & \OW.
\ecases
\]
We show that $w$ is a BJ supersolution of
\beq\label{hje}
u_t+H(x,D_xu,u)=0 \ \ \IN M\tim(0,\infty).
\eeq

To do this, we note that the function $u(x):=d(x,y)$ is a solution of the eikonal
equation
\[
|Du(x)|=1 \quad \IN M\stm\{y\}.
\]
This is a standard observation, and skip the proof here.  From this remark,
we find that the function $u(x,t):=d(x,y)-\gd t$  is a solution of
\beq \label{propagation}
u_t+\gd|D_xu|=0 \ \ \ \IN (M\stm\{y\})\tim(0,\infty).
\eeq
Now, we fix $\gth\in C^1(\R,\R)$ such that
\[
\gth(r)=0 \ \ \FOR\  r\leq 0,\quad
\gth(r)>0 \ \ \FOR\ r>0, \ \ \AND \ \
\gth'(r)\geq 0 \ \ \FOR\  r\in\R,
\]
define $z^k\in C(M\tim[0,\infty),\R)$ for every $k\in\N$ by
\[
z^k(x,t)=k\gth(d(x,y)-\gd t),
\]
and observe that for any $k\in\N$, $z^k$ is a solution of \erf{propagation}.
Moreover, by sending $k\to\infty$, we find that the function
$z\in \LSC(M\tim[0,\infty),\R\cup\{\infty\})$ given by
\[
z(x,t)=
\bcases
0 & \IF \ d(x,y)\leq\gd t,\\
+\infty & \OW,
\ecases
\]
is a solution (in the BJ sense) of \erf{propagation}.  Noting that the set
$Z:=\{(x,t)\in M\tim[0,\infty)\mid d(x,y)\leq \gd t\}$ is a {neighborhood}
of $\{y\}\tim(0,\infty)$ and $z$ vanishes on $Z$, we easily deduce  that
$z$ is a {solution} of
\[
u_t+\gd|D_xu|=0 \ \ \IN M\tim(0,\infty).
\]
It is now easy to check that the function $w(x,t)=z(x,t)+\gL^{-1}(e^{\gL t}-1)$
is a solution of
\[
u_t+\gd |D_x u|-C_1-\gL|u|=0 \ \ \IN M\tim(0,\infty),
\]
which assures, together with \erf{linear-growth}, that $w$ is a supersolution
of \erf{hje}. Since $w(x,0)=z(x,0)=h(x,0,y,0)$ for all $x\in M$, we conclude by Theorem \ref{cwch} that $h(x,t,y,0)\leq w(x,t)$ for all $(x,t)\in M\tim[0,\infty)$, which yields \erf{u-bound}.
\eproof

\section{Long-time behavior of solutions}

We are concerned with the long-time behavior of the solution $u=u(x,t)$ of problem
\beq\label{hje-lt}
u_t+H(x,D_xu,u)=0 \ \ \IN \ M\tim(0,\infty).
\eeq

\begin{theorem}\label{finiteness} Assume \erf{H1}--\erf{H4}. Let $u\in \LSC(M\tim(0,\infty),\R\cup\{\infty\})$
be a BJ solution of \erf{hje-lt}. Set
\beq \label{def-v}
v(x)=\lim_{r\to 0+}\inf\{u(y,t)\mid d(y,x)<r,\ t>r^{-1}\} \ \ \FOR\ x\in M.
\eeq
Assume that $v(z)\in\R$ for some $z\in M$. Then, $v\in \Lip(M)$ and $v$ is a
viscosity solution of $H(x,Dv,v)=0$ in $M$.
\end{theorem}

We remark, as noted in the proof of Lemma \ref{l1.cwch}, that, once the Lipschitz continuity of $v$ is known,
$v$ is a BJ solution of  $H(x,Dv,v)=0$ in $M$ if and only if it is a
viscosity solution, in the Crandall-Lions sense, of $H(x,Dv,v)=0$ in $M$.

Notice that the definition of $v$ above is the so-called lower relaxed limit
of $u(x,t)$ as $t\to \infty$. In particular, $v$ is lower semicontinuous on $M$.

An immediate consequence of the theorem above is the following.

\begin{corollary} \label{cor.finite} Assume \erf{H1}--\erf{H4}. Let
$v\in \BJ(H,M)$ If $v(z)\in\R$ for some $z\in M$, then $v\in\Lip(M,\R)$.
\end{corollary}

Notice that, thanks to the above corollary, if $v\in\BJ(H,M)$ and $v(x)\not\equiv +\infty$, then $v\in\Lip(M,\R)$ and, consequently, $v\in\cS(H,M)$.

\bproof Set $u(x,t)=v(x)$ for $(x,t)\in M\tim(0,\infty)$ and note that
$u\in \BJ(\pl_t-H, M\tim(0,\infty))$. For these functions $u$ and $v$, the relation
\erf{def-v} holds and, hence, Theorem \ref{finiteness} assures that $v\in\Lip(M,\R)$.
\eproof

\bproof[Proof of Theorem \ref{finiteness}]  We first show that $v(x)\in\R$ for all $x\in M$. For this, we set
$\gam:=\inf_{M}v<\infty$ and pick
a sequence $(y_k,t_k)\in M\tim(0,\infty)$ such that $\lim_{k\to\infty}t_k=\infty$
and $\lim_{k\to\infty} u(y_k,t_k)=\gam$.
 According to Proposition \ref{two-bounds}, there exist constants $\gd>0,\, C>0$
such that for all $x,y\in M,\, t\geq 0$,
\[
h(x,t,y,0)\leq C(e^{\gL t}-1) \ \ \IF \ d(x,y)\leq \gd t.
\]
We fix $T>0$ so that $\gd T$ is larger than or equal to the diameter of $M$, that is,
\[
d(x,y)\leq \gd T \ \ \FORALL x,y\in M.
\]
Hence, we have
\beq\label{bound-u+}
h(x,T,y,0)\leq C(e^{\gL T}-1) \ \ \FORALL x,y\in M.
\eeq
By Theorem \ref{f2.thm}, we have
\beq \label{f-formula}
u(x,t_k+T)=(S_T u(\cdot,t_k))(x) \leq h(x,T,y_k,u(y_k,t_k)) \ \ \FORALL x\in M.
\eeq

Now, we suppose that $\gam=-\infty$. We may assume without loss of generality
that $u(y_k,t_k)\leq -k$ for all $k\in\N$. Combining this with \erf{f-formula},
Lemma \ref{f1.lemma}, Corollary \ref{cor-cwch}, and \erf{bound-u+},
we obtain for all $x\in M$ and $k\in \N$,
\[
u(x,t_k+T)\leq h(x,T,y_k,-k)\leq -k e^{-\gL T}+h(x,T,y_k,0)
\leq -k e^{-\gL T}+C(e^{\gL T}-1).
\]
This shows that $v(x)=-\infty$ for all $x\in M$, which contradicts that $v(z)\in\R$.
Hence, we find that $\gam\in \R$.

We choose a constant $C_0>0$ so that $u(y_k,t_k)\leq C_0$ for all $k\in\N$.
We argue similarly to the above by
using \erf{f-formula}, Lemma \ref{f1.lemma}, Corollary \ref{cor-cwch}, and
\erf{bound-u+}, to obtain for all $x\in M$ and $k\in\N$,
\[
u(x,t_k+T)\leq h(x,T,y_k,C_0)\leq C_0 e^{\gL T}+C(e^{\gL T}-1),
\]
which implies
\beq \label{bound-u++}
v(x)\leq C_0 e^{\gL T}+C(e^{\gL T}-1) \ \ \FORALL\ x\in M.
\eeq
Thus, we conclude that $v(x)\in\R$.

As a basic property of the lower half-limit,
recalling the definition of $v$ and noting that the functions $(x,t) \mapsto u(x,t+s)$, with $s>0$, are BJ solutions
of  \erf{hje-lt},
we find that
$v$ is a viscosity supersolution, in the Crandall-Lions sense,
both of the Hamilton-Jacobi equations $\pl_t u+H(x,D_xu,u)=0$ and $-\pl_t u-H(x,D_xu,u)=0$ in $M\tim(0,\infty)$, which means that
$v$ is a BJ solution of $H(x,Du,u)=0$ in $M$.

Since $v\in\LSC(M,\R)$, it follows that $v$ is bounded from below on $M$,
which, together with \erf{bound-u++}, assures that $v$ is bounded on $M$.
Let $C_1>0$ be a constant such that $|v(x)|\leq C_1$ for all $x\in M$.
If we set $w=-v\in \USC(M,\R)$, then $w$ is a viscosity subsolution, in the Crandall-Lions sense, of $H(x,-Du,-u)=0$ in $M$. Since
$H(x,p,u)\geq H(x,p,0)-\gL C_1$ as far as $|u|\leq C_1$, the Hamiltonian
$(x,p)\mapsto H(x,p,0)-\gL C_1$ is coercive and $v$ is a subsolution, in the Crandall-Lions sense, of
$H(x,-Du,0)-\gL C_1=0$ in $M$, a standard regularity result assures that $v\in\Lip(M,\R)$.

The proof is now complete.
 \eproof

Given a Hamiltonian $H\in C(T^*M\tim\R,\R)$, we set
\[
H^\ominus(x,p,u)=H(x,-p,-u).
\]
Note that $(H^\ominus)^\ominus=H$ and
that if $H$ satisfies \erf{H1}--\erf{H4}, then so does $H^\ominus$.

Under the hypotheses \erf{H1}--\erf{H4}, we write $S_t$ and $S_t^\ominus$ for
the operators $S_t^H$ and $S_t^{H^\ominus}$, respectively.

\begin{remark}
By \cite[Proposition 2.8]{WWY1}, for any $\phi\in C(M,\R)$, there hold
\[S_t\phi=T_t^-\phi,\quad S_t^\ominus\phi=-T_t^+(-\phi),\]
where $T_t^{-}$ and $T_t^+$ denote the backward and forward Lax-Olenik semigroup associated with $H$, respectively. By using  $T_t^{\pm}$, new progress on viscosity solutions of contact Hamilton-Jacobi equations was achieved (\cite{WWY1a,WWY1b,NWY2021}). We also refer the reader to \cite{WY2021} for existence and long time behavior of viscosity solutions of contact Hamilton-Jacobi  equation where no monotonicity assumptions is imposed. Besides,  by using the Herglotz variational principle (\cite{CCJWY}), some kinds of representation formulae for the viscosity solution of (\ref{hje-lt}) on the Cauchy problem were  provided in  \cite{HCHZ}.
\end{remark}

We establish the following theorems.

\begin{theorem}\label{thm1} Assume \erf{H1}--\erf{H4}. For any $u_0\in\cS(H)$,
the function $t\mapsto S_t^\ominus (-u_0)(x)$ is nondecreasing on $[0,\infty)$
for any $x\in M$, and the limit
\[
v_0(x):=\lim_{t\to\infty} S_t(-u_0)(x)
\]
exists for any $x\in M$, and $v_0\in \cS(H^\ominus)$.
The convergence above is uniform on $M$.
\end{theorem}

In view of the theorem above, under the hypotheses (H1)--(H4), we may introduce
the operators $T_\infty\mid \cS(H^\ominus)\to\cS(H),\,T_\infty^\ominus\mid
\cS(H)\to\cS(H^\ominus)$ by
\[\bald
T_\infty v(x)&=\lim_{t\to\infty}S_t(-v)(x) \ \ \FOR x\in M \AND v\in \cS(H^\ominus),
\\T_\infty^\ominus u(x)&=\lim_{t\to\infty}S_t^\ominus(-u)(x) \ \ \FOR x\in M \AND u\in \cS(H).
\eald
\]
The monotonicity of $S_t(-v)(x)$ and $S_t^\ominus(-u)(x)$ in $t$ yields
\beq\label{mono} \left\{\bald
&T_\infty v\geq -v \ \ \ON M \ \FORALL v\in \cS(H^\ominus),
 \\&T_\infty^\ominus u\geq -u \ \ \ON M \ \FORALL u\in\cS(H).
\eald\right. \eeq
Also, the comparison principle implies
\beq\label{compT}\left\{
\bald
& v_1,v_2\in \cS(H^\ominus),\ v_1\leq v_2 \ \implies \ T_\infty v_1\geq T_\infty v_2,
\\& u_1,u_2\in\cS(H),\ u_1\leq u_2 \ \implies \ T_\infty^\ominus u_1\geq T_\infty^\ominus u_2.
\eald\right.
\eeq

Let $\cI(T_\infty)$ and $\cI(T_\infty^\ominus)$ denote the images of
$T_\infty$ and $T_\infty^\ominus$, respectively, that is,
\[
\cI(T_\infty)=\{T_\infty v\mid v\in\cS(H^\ominus)\}
\ \ \AND \ \ \cI(T_\infty^\ominus)=\{T_\infty^\ominus u \mid u\in\cS(H)\}.
\]

\begin{theorem} \label{thm2} Assume \erf{H1}--\erf{H4}. \begin{enumerate} \item[(1)] $\cS(H)
\not=\emptyset$\ if and only if\ $\cS(H^\ominus)\not=\emptyset$.

\item[(2)] For any $u\in \cI(T_\infty)$ and $v\in\cI(T_\infty^\ominus)$,
\[
T_\infty\circ T_\infty^\ominus u=u \ \ \AND \ \ T_\infty^\ominus\circ
T_\infty v=v.
\]

\item[(3)] Let $u\in\cI(T_\infty)$ and $\phi\in\LSC(M,\R)$. Assume that there exist a finite number of $v_1,\ldots, v_k \in\cS(H^\ominus)$ such that \
$T_\infty v_i=u$ \ for all $i=1,\ldots,k$,
and  $\min_i (-v_i)\leq \phi \leq u$ on $M$. Then
\[
\lim_{t\to \infty} S_t \phi=u\ \ \text{ uniformly on }M.
\]
\end{enumerate}
\end{theorem}

For the proof of the theorems above, we need some lemmas.
In the following lemmas, we always \emph{assume} \erf{H1}--\erf{H4}.

\begin{lemma} \label{extremal} Let $u_0\in \cS(H)$ and $x\in M$. There exists a curve $\gam : (-\infty, 0]\to M$ such that $\gam(0)=x$, $\gam\in \AC([-\tau,0],M)$
for every $\tau>0$, and, for all $t>0$,
\[
u_0(\gam(0))=u_0(\gam(-t))+\int_{-t}^0 L(\gam(s),\dot\gam(s),u_0(\gam(s)))ds.
\]
\end{lemma}

The lemma above is a classical observation and follows from Theorem \ref{thm:repr}.

\bproof[Outline of proof] Since $u_0\in\BJ(\pl_t+H)$, by Theorem \ref{thm:repr}, there exists
$\{\gam_k\}_{k\in\N}\subset \AC([-1,0],M)$ such that
\[
u_0(\gam_k(0))=\int_{-1}^0 L(\gam_k(s),\dot\gam_k(s),u_0(\gam_k(s))ds+u_0(\gam_k(-1)) \ \ \FORALL \ k\in\N,
\]
\[
\gam_1(0)=x, \quad \AND\quad \gam_{k+1}(0)=\gam_k(-1) \ \ \FORALL k\in\N.
\]
Define $\gam : (-\infty,0] \to M$ by
\[
\gam(s)=
\bcases
\gam_1(s) & \FOR s\in(-1,0], \\
\gam_2(s+1)&\FOR s\in(-2,-1],\\
\gam_3(s+2)&\FOR s\in(-3,-2],\\
 \ \ \ \  \vdots & \quad \ \ \vdots
\ecases
\]
It is easily checked that $\gam$ has all the required properties in Lemma \ref{extremal}.
\eproof

\begin{lemma} \label{estimate}Let $u_0\in\cS(H)$ and set $v_0=-u_0$.
Then, $v_0\in \BJ^-(H^\ominus)$ and
\[
\min_{x\in M} (S_tv_0(x)-v_0(x))=0 \ \ \FORALL t>0.
\]
\end{lemma}

\bproof We set $v(x,t)=S_t v_0(x)$ for $(x,t)\in M\tim[0,\infty)$.

Since $u_0\in\Lip(M,\R)\cap\cS^-(H)$,
we have
\[
v_0\in \BJ^-(H^\ominus)\subset \BJ^-(\pl_t+H^\ominus),
\]
and, by Theorem \ref{cwch},
\beq\label{mono1}
v_0(x)\leq v(x,t) \ \ \FORALL (x,t)\in M\tim[0,\infty).
\eeq

Fix any $x$, and, in view of Lemma \ref{extremal}, choose a curve $\gam$, with $\gam(0)=x$, such that for all $t>0$,
\[
\gam\in\AC([-t,0],M) \ \ \AND \ \
u_0(\gam(0))=u_0(\gam(-t))+\int_{-t}^0 L(\gam(s),\dot\gam(s),u_0(\gam(s)))ds.
\]
Set $\eta(s):=\gam(-s)$, to find that for all $t>0$,
\[
u_0(\eta(0))=u_0(\eta(t))+\int_0^t L(\eta(s),-\dot\eta(s),u_0(\eta(s)))ds,
\]
which reads
\[
v_0(\eta(t))=v_0(\eta(0))+\int_0^t L(\eta(s),-\dot\eta(s),-v_0(\eta(s)))ds.
\]

By Theorem \ref{thm:repr}, we have
\[
v(x,t)=\inf_{\xi(t)=x}\left(v_0(\xi(0))+\int_0^t L(\xi(s),-\dot\xi(s),-v(\xi(s),s))ds\right).
\]
In particular,
\[
v(\eta(t),t)\leq v_0(\eta(0))+\int_0^t L\left(\eta(s),-\dot\eta(s),-v(\eta(s),s))\right)ds.
\]
Since $u\mapsto L(x,\xi,-u)$ is Lipschitz continuous in $\R$, with $\gL$ as a Lipschitz bound, and, by \erf{mono1}, $v(\eta(s),s)\geq v_0(\eta(s))$ for all $s\geq 0$, we have
\[\bald
L(\eta(s),-\dot\eta(s),-v(\eta(s),s))
\leq L(\eta(s),-\dot\eta(s),-v_0(\eta(s)))+\gL(v(\eta(s),s)-v_0(\eta(s)).
\eald
\]
Thus, for all $t>0$,
\[
v(\eta(t),t)\leq v_0(\eta(t))+\int_0^t\gL(v(\eta(s),s)-v_0(\eta(s)))ds,
\]
which implies that
\[
v(\eta(t),t)-v_0(\eta(t))\leq 0 \ \ \FORALL t>0.
\]
We conclude that
\[
\min_{x\in M}(v(x,t)-v_0(x))=0 \ \ \FORALL t\geq 0. \qedhere
\]
\eproof

\begin{lemma} \label{nondec}  For any $\phi\in \BJ^-(H)$,  the function $t\mapsto S_t \phi (x)$ is nondecreasing
on $[0,\infty)$ for every $x\in M$.
\end{lemma}

\bproof Fix any $\phi\in \BJ^-(H)$ and set $u(x,t)=S_t\phi(x)$ for
$(x,t)\in M\tim[0,\infty)$.

By Theorem \ref{cwch}, we have
\[
\phi(x)\leq u(x,t)  \ \ \FORALL (x,t)\in M\tim[0,\infty).
\]
For any $\gd>0$,
\[
\phi(x)\leq u(x,\gd) \ \ \FORALL x\in M,
\]
and hence, again by the comparison priciple,
\[
u(x,t)\leq (S_t u(\cdot,\gd))(x)=u(x,\gd+t) \ \ \FORALL (x,t)\in M\tim [0,\infty).
\]
This shows that the function $t\mapsto S_t\phi(x)$ is nondecreasing on $[0,\infty)$
for any $x\in M$.
\eproof

\bproof[Proof of Theorem \ref{thm1}] Fix any $u_0\in\cS(H)$ and set
$\phi=-u_0$.
Set
\[
v_0(x):=\lim_{r\to 0+}\inf\{S_t^\ominus \phi (y)\mid d(y,x)<r,\ t>r^{-1}\}
 \ \ \FOR x \in M.
\]
By Lemma \ref{estimate}, we have $\phi\in\BJ^-(H^\ominus)$ and
\[
\min_{x\in M}(S_t^\ominus \phi(x)-\phi(x))=0 \ \ \FORALL t>0,
\]
which implies
that, since $M$ is compact,
\[
v_0(z)<\infty \ \ \text{ for some }z\in M.
\]
Theorem \ref{finiteness}, with $H$ replaced by $H^\ominus$,
ensures that $v_0\in \Lip(M,\R)$ and $v_0\in\cS(H^\ominus)$.

By Lemma \ref{nondec}, with $H^\ominus$ in place of $H$,
the function $t\mapsto S_t^\ominus \phi(x)$ is nondecreasing on $[0,\infty)$
for any $x\in M$. Hence, $S_\tau^\ominus\phi\leq S^\ominus_{\tau+t} \phi$ on $M$ for all $\tau,t>0$, which
shows that $S_\tau^\ominus \phi\leq v_0$ on $M$ for all $\tau>0$.
On the other hand, by the definition of $v_0$, we have \
$\lim_{t\to\infty}S_t^\ominus\phi(x)\geq v_0(x)$ for all \ $x\in M$.
Thus, we find that
\[
\lim_{\tau\to \infty}S_\tau^\ominus \phi(x)=v_0(x) \ \ \FORALL\ x\in M.
\]
The convergence above is uniform in $x\in M$ by the Dini theorem.
\eproof

\bproof[Proof of Theorem \ref{thm2}] We first treat (1).  It is an immediate consequence
of Theorem \ref{thm1} that $\cS(H)\not=\emptyset$ implies $\cS(H^\ominus)\not= \emptyset$. The converse implication also follows from Theorem \ref{thm1}, with
$H$ replaced by $H^\ominus$.

We next consider (2).  Fix any  $u\in\cI(T_\infty)$ and choose
$v_0\in \cS(H^\ominus)$ such that
$T_\infty v_0=u$. Set $v:=T_\infty^\ominus u$.
It follows from \erf{mono} that $v_0+u\geq 0$ and $u+v\geq 0$ on $M$.
Hence, by the comparison principle, we find that
$S^\ominus_t(-u)\leq v_0$  and $S_t(-v)\leq u$ on $M$ for all $t\geq 0$, which implies that $T_\infty^\ominus u\leq v_0$ and $T_\infty v\leq u$ on $M$.
From the first of these inequalities, together with \erf{compT}, we obtain $T_\infty\circ T_\infty^\ominus u\geq T_\infty v_0=u$ on $M$, which, combined with the latter of the inequalities above, yields
$u\geq T_\infty v=T_\infty\circ T_\infty^\ominus u\geq T_\infty v_0=u$ on $M$.
This shows that $T_\infty\circ T_\infty^\ominus u=u$ on $M$.
Replacing $H$ by $H^\ominus$ and noting that $(H^{\ominus})^\ominus$, we also
conclude that $T_\infty^\ominus\circ T_\infty v=v$ on $M$ for all
$v\in\cI(T_\infty^\ominus)$.

Now, we treat  (3).  First of all, we show that for all $(x,t)\in M\tim[0,\infty)$,
\beq \label{pr-(3)}
\min_{1\leq i\leq k}S_t(-v_i)(x)=S_t(\min_{1\leq i\leq k}(-v_i))(x).
\eeq

Indeed, by Lemma \ref{BJ-min}, we have
\[
(x,t)\mapsto \min_{1\leq i\leq k}S_t(-v_i)(x) \, \in \, \BJ(\pl_t+H,M\tim(0,\infty)).
\]
Moreover, it is easily checked that
\[
(x,t)\mapsto \min_{1\leq i\leq k}S_t(-v_i)(x) \in\LSC(M\tim[0,\infty),\R\cup\{\infty\}),
\]
and for all $x\in M$,
\[
\min_{i\leq i\leq k}(-v_i)(x)=\liminf_{t\to \infty}\min_{1\leq i\leq k}S_t(-v_i)(x).
\]
Thus, by the definition of the operator $S_t$, we conclude that \erf{pr-(3)} holds.

By the choice of $v_i$, we have
\[
\lim_{t\to\infty}S_t(-v_i)=u \ \ \text{ uniformly on }M \ \FORALL i\in\{1,\ldots,k\},
\]
which assures together with \erf{pr-(3)} that
\beq\label{pr-(3)+}
\lim_{t\to\infty} S_t(\min_{1\leq i\leq k}(-v_i))(x)=u(x)  \ \ \text{ uniformly on }M.
\eeq
Since $\min_{1\leq i\leq k}(-v_i)\leq\phi\leq u$ on $M$ and $u\in\BJ(H)$, we deduce by the comparison principle that for all $t\geq 0$,
\[
S_t(\min_{1\leq i\leq k}(-v_i))\leq S_t\phi\leq S_t u=u \ \ \ON M,
\]
which, combined with \erf{pr-(3)+}, yields
\[
\lim_{t\to\infty} S_t\phi =u  \ \ \text{ uniformly on } M.  \qedhere
\]
\eproof

\section{Classification of the solutions of $H=0$}

For $v\in \cS(H^\ominus)$, we define the set $\mathcal{D}_v$ by
\[\mathcal{D}_v:=\{w\in \cS(H)\mid T_\infty^{\ominus}w=v\}.\]

The sets $\cD_v$, with $v\in\cI(T_\infty^\ominus)$, classifies $\cS(H)$ as follows.

\begin{theorem} \label{classify} Assume \erf{H1}--\erf{H4}. Then:
\begin{enumerate}\item[(1)] $\#\cI(T_\infty)=\#\cI(T_\infty^\ominus)$,
where $\# A$ denotes the cardinality of set $A$.
\item[(2)] ${\displaystyle \cS(H)=\bigsqcup_{v\in \cI(T_\infty^\ominus)}\mathcal{D}_v
=\bigsqcup_{u\in \cI(T_\infty)} \cD_{T_\infty^\ominus u}}$.
\item[(3)] $w\geq u$ for all $u\in\cI(T_\infty)$ and $w\in \cD_{T_\infty^\ominus u}$.
\end{enumerate}
\end{theorem}

\bproof  By (2) of Theorem \ref{thm2}, $T_\infty^\ominus$ is a bijection of $\cI(T_\infty)$ to $\cI(T_\infty^\ominus)$, with the inverse map $T_\infty$. Hence,
$\#\cI(T_\infty^\ominus)=\#\cI(T_\infty)$, which proves (1).

Since $T_\infty^\ominus\mid \cS(H)\to \cI(T_\infty^\ominus)$
is a surjection,  if we introduce the binary relation $\sim$ on $\cS(H)$ by
\[
u_1\sim u_2 \ \iff \ T_\infty^\ominus u_1=T_\infty^\ominus u_2,
\]
then this relation is clearly an equivalence relation on $\cS(H)$ and the sets
$\mathcal{D}_v$, with
$v\in\cI(T_\infty^\ominus)$, constitute all the equivalence classes in this relation. Consequently,
\[
\cS(H)=\bigsqcup_{v\in \cI(T_\infty^\ominus)}\mathcal{D}_v.
\]
Since $T_\infty\mid \cI(T_\infty^\ominus) \to \cI(T_\infty)$ is bijective, we have
\[
\cS(H)=\bigsqcup_{v\in \cI(T_\infty^\ominus)}\mathcal{D}_v
=\bigsqcup_{u\in \cI(T_\infty)}\mathcal{D}_{T_\infty^\ominus u},
\]
which proves assertion (2).

To check (3), let $u\in\cI(T_\infty)$ and $w\in \cD_{T_\infty^\ominus u}$. By \erf{mono}, we have $T_\infty^\ominus w+w\geq 0$ on $M$ and, hence, $w=S_t w\geq S_t(-T_\infty^\ominus w)=S_t(-T_\infty^\ominus u)$ on $M$ for all $t\geq 0$, which implies that $w\geq T_\infty\circ T_\infty^\ominus u=u$ on $M$.
\eproof

\begin{remark} (1) Under the assumptions \erf{H1}--\erf{H4} and that $\cS(H)\not=\emptyset$, if $H=H(x,p,u)$ is strictly monotone in $u$, then $\#\cI(T_\infty)=1$. More precisely,
\begin{itemize}
\item [(i)] if $H$ is strictly increasing in $u$, then
\[\#\cI(T_\infty)=\#\cS(H)=1;\]
\item [(ii)] if $H$ is strictly decreasing in $u$, then
\[\#\cI(T_\infty)=1\leq \#\cS(H).\]
\end{itemize}
Concerning (ii) above, consider two examples. The first example is about the equation
(see \cite[Proposition 10]{DW})
\begin{equation}%\tag{E1}
\label{e1}
 -u+\frac{1}{2}|Du|^2+\cos x-1=0\quad \IN \mathbb{T}.
\end{equation}
For \erf{e1}, we have $\#\cS(H)=1.$
The other one concerns the equation (see \cite[Example 1.1]{WWY1})
\begin{equation}%\tag{E2}
\label{e2}
-u+\frac{1}{2}|Du|^2=0\quad \IN \mathbb{T}.
\end{equation}
For \erf{e2}, the solutions $u$ are given by
\[
u(x)=\min_{y\in K} \frac 12 d(x,y)^2,
\]
with $K$ being nonempty compact subsets of $\mathbb{T}$.
Since the totality of compact subsets of $\T$ is an infinite set,
we have $\#\cS(H)=\infty. $
Moreover, one can choose a Hamiltonian $H$ so that
$\#\cI(T_\infty)<\#\cS(H)<\infty$.
The following equation is taken from \cite[Proposition 1.14]{WWY2}:
\begin{equation}%\tag{E3}
\label{e3}
 -u+\frac{1}{2}|Du|^2+Du\cdot V(x)=0\quad \IN \mathbb{T},
\end{equation}
where $V:\mathbb{T}\to \R$  is a smooth function which has exactly two zeros $x_1$, $x_2$ with $V'(x_1)>0$, $V'(x_2)<0$. For (\ref{e3}), we have
\[\#\cS(H)=2.\]
\noindent {(2)}
 If $H$ is non-monotone in $u$, then it may happen that $\#\cI(T_\infty)>1$. For example,
consider the equation \[f(u)+\frac{1}{2}|Du|^2=0\quad \IN \mathbb{T},\]
where $f$ is a smooth function  satisfying $f(u)=f(-u)$,  $f> 0$ on $(-1,1)$ and $f(u)=u+1$ for $u\in (-\infty, -\frac{1}{2}]$.
In this case, the solutions $u$ are given by either $u(x)=-1$, or
\[
u(x)=1+\min_{y\in K}\frac 12 d(x,y)^2,
\]
where $K$ ranges over all compact subsets of $\T$, and we have
\[\#\cI(T_\infty)=2, \quad \#\cS(H)=\infty.\]
The Hamiltonian $H$ associated with \erf{e3} is self-adjoint in the sense that
$H^\ominus=H$.
\noindent{(3)} If $H$ is independent of $u$ and $\cS(H)\not=\emptyset$, then
\[\#\cI(T_\infty)=\#\cS(H)=\infty.\]
Moreover, the structure of $\cS(H)$ can be described in terms of static classes in the Aubry set (see \cite[Theorem 0.2]{Co}).
\end{remark}

\appendix
\section{A classical existence result}

We give here a proof of the following classical existence theorem for a viscosity solution
of \erf{cp.1}.

\begin{theorem} \label{a.1} Assume \erf{H1}--\erf{H4} and that $\phi\in {C^1(M,\R)}$.
Then there exists a unique viscosity solution $u\in\Lip(M\tim[0,T),\R)$ of \erf{cp.1}.
\end{theorem}

We begin with a lemma concerning the Cauchy problem for
\beq\label{a.2}
u_t+H(x,D_xu,u)+\gam|D_xu|^2=0 \ \ \IN M\tim(0,T),
\eeq
where $\gam$ and $T$ are constants such that $0\leq \gam<\infty$
and $0<T\leq\infty$.

\begin{lemma} \label{a.3} Assume that $H$ is a bounded and uniformly continuous
function on {$T^*M \tim\R$} and satisfies \erf{H3} and \erf{H4}.
Let $v\in\USC(M\tim[0,T),\R)$ and $w\in \LSC(M\tim[0,T),\R)$ be a viscosity subsolution and supersolution, respectively, of \erf{a.2}.
Assume that $v,\,w$ are locally bounded on $M\tim[0,T)$ and that $v(x,0)\leq w(x,0)$
for all $x\in M$. Then $v\leq w$ on $M\tim[0,T)$.
\end{lemma}

The proof below is similar to but simpler than that for Theorem \ref{cwch}.

\bproof  We may assume without loss of generality that $T<\infty$.
If we set
\[
\tilde v(x,t)=e^{-\gL t}v(x,t) \ \ \AND \ \ \tilde w(x,t)=e^{-\gL t}w(x,t)
\ \ \FOR (x,t)\in M\tim [0,T),
\]
then $\tilde v$ and $\tilde w$ are a viscosity subsolution and supersolution,
respectively, of
\[
u_t+\hH(x,D_xu, t,u)+\gam e^{\gL t}|D_xu|^2=0,
\]
where $\tH\in C(T^*M\tim[0,T)\tim\R,\R)$ is given by
\[
\hH(x,p,t,u):=\gL u+e^{-\gL t}H(x,e^{\gL t}D_xu, e^{\gL t}u).
\]
We denote by $\tH$ the function $(x,p,t,u)\mapsto e^{-\gL t}H(x,e^{\gL t}p,e^{\gL t}u)$ on $T^*M\tim[0,T)\tim \R$.
Note that for any $(x,p,t)\in T^*M\tim[0,T)$, the function $u\mapsto
\hH(x,p,t,u)$ is nondecreasing on $\R$ and the function $u\mapsto
\tH(x,p,t,u)$ is Lipschitz continuous with Lipschitz bound $\gL$.

We need only to show that $\tilde v\leq \tilde w$ on $M\tim[0,T)$, and
by contradiction, we suppose that $\sup_{M\tim[0,T)}(v-w)>0$.
For ease of notation, we henceforth write $v$ and $w$ for $\tilde v$ and $\tilde w$, respectively.

For $\ep>0$, set
\[
v^\ep(x,t)=v(x,t)-\fr{\ep }{T-t+\ep^2}  \ \ \FOR \ (x,t)\in M\tim[0 T).
\]
Note that, as $\ep\to 0+$, $\sup_{M\tim (T-\ep^2,\,T)}v^\ep(x,t) \to -\infty$
and $v^\ep(x,t) \to v(x,t)$ uniformly on $M\tim [0,S]$ for all $0<S<T$, and
that $v^\ep(x,0)<w(x,0)$ for all $x\in M$. We may fix $\ep>0$ small enough so that $v_\ep-w$ takes a positive maximum at
some point $(x_0,t_0)\in M\tim(0,T)$.
It is easily seen that $v^\ep$ is a viscosity subsolution of
\[
u_t+\hH(x,D_xu,t,u)+\gam e^{\gL t}|D_xu|^2+\fr{\ep}{(T-t+\ep^2)^2}=0 \ \ \IN M\tim(0,T).
\]
Setting \ $\gd=\ep/(2(T+\ep^2)^2)$, we see immediately that $v^\ep$
is a viscosity subsolution of
\beq\label{a.3.0}
u_t+\hH(x,D_xu,t, u)+\gam e^{\gL t}|D_xu|^2+2\gd=0 \ \ \IN M\tim(0,T).
\eeq

We fix a local coordinates in an open neighborhood $V$ of the point $x_0$
so that we may regard $V$ as an open subset of $\R^n$
and $T^*V=V\tim\R^n$.  We introduce a (bump) function
$\rho\in C^1(V\tim(0,T),\R)$ such that
\beq\label{a.3.1}\left\{\ \bald
&(\rho-w)(x_0,t_0)>0,
\\&(\rho-w)(x,t)<-1 \ \ \FORALL \ (x,t)\in V\tim(0,T)\stm B,
\eald \right.
\eeq
where $B$ is a bounded
open subset of $\R^{n+1}$ whose closure $\ol B$ is included in $V\tim(0,T)$.

For $\eta\in (0,1)$, we set
\[
v^{\ep,\eta}(x,t)=(1-\eta)v^\ep(x,t)+\eta\rho(x,t) \ \ \FOR \ (x,t)\in V\tim(0,T),
\]
and note by \erf{a.3.1} that the function
\[
v^{\ep,\eta}(x,t)-w(x,t)=(1-\eta)(v^\ep-w)(x,t)+\eta(\rho-w)(x,t)
\]
takes a positive maximum at some point $(x_\eta,t_\eta)\in B$ and
\beq\label{a.3.2}
v^{\ep,\eta}(x,t)-w(x,t)<v^{\ep,\eta}(x_\eta,t_\eta)-w(x_\eta,t_\eta)-\eta
\ \ \FORALL \ (x,t)\in V\tim(0,T)\stm B.
\eeq

Taking into account the convexity of $p\mapsto
\tH(x,p,u)+\gam e^{\gL t}|p|^2$ and the Lipschitz property
of $\tH$, we compute in a slightly formal way that
\[\bald
v_t^{\ep,\eta}&+\hH(x,D_x v^{\ep,\eta},t,v^{\ep,\eta})
+\gam e^{\gL t}|D_xv^{\ep,\eta}|^2
\\&\leq (1-\eta)(v^\ep_t+\gL v^\ep+\tH(x, D_xv^\ep,t,v^{\ep,\eta})
+\gam e^{\gL t}|D_xv^{\ep}|^2)
\\&\quad +\eta(\rho_t+\gL \rho+\tH(x,D_x\rho,t,v^{\ep,\eta})
+\gam e^{\gL t}|D_x\rho|^2)
\\&\leq (1-\eta)(v^\ep_t+\gL v^\ep+\tH(x, D_xv^\ep,t,v^\ep)
+\gam e^{\gL t}|D_xv^{\ep}|^2
+\eta\gL(|v^\ep|+|\rho|)))
\\&\quad +\eta(\rho_t+H(x,D_x\rho,\rho)+\gam|D_x\rho|^2+\gL(|v^\ep|+|\rho|)).
\eald
\]
Hence,
thanks to \erf{a.3.0}, we may choose $\eta\in(0,1)$ small enough so that
$v^{\ep,\eta}$ is a viscosity subsolution of
\[
u_t+\hH(x,D_xu,t,u)+\gam e^{\gL t}|D_xu|^2+(1-\eta)\gd=0 \ \ \IN V\tim (0,T).
\]

Consider the function
\[
v^{\ep,\eta}(x,t)-w(y,s)-\fr 12\left(\ga|x-y|^2+\ga^2(t-s)^2\right)
\]
on $\ol B\tim\ol B$, where $\ga>1$, and pick a maximum point
$(\hat x,\hat t,\hat y,\hat s)$ of this function. It is a standard observation that, as
$\ga\to \infty$,
\beq\label{a.3.3}
\ga|\hat x-\hat y|^2+\ga^2(\hat t-\hat s)^2 \to 0,
\eeq
and that for any limiting point $(\bar x,\bar t,\bar y,\bar s,\bar v,\bar w)$, as
$\ga \to \infty$,  of the family
\[
\{(\hat x,\hat t,\hat y,\hat s, v^{\ep,\eta}(\hat x,\hat t),w(\hat y,\hat s))\},
\]
we have
\beq\label{a.3.4}
(\bar x,\bar t)=(\bar y,\bar s) \ \ \AND \ \ \bar v-\bar w=\max_{\ol B}(v^{\ep,\eta}-w).
\eeq
Fix such a limiting point $(\bar x,\bar t,\bar y,\bar s,\bar v,\bar w)$ and a
sequence $\{\ga_j\}_{j\in\N}$ such that, as $j\to\infty$, the sequence of
the points
\[
(\hat x,\hat t,\hat y,\hat s, v^{\ep,\eta}(\hat x,\hat t),w(\hat y,\hat s)),
\]
with $\ga=\ga_j$, converges to $(\bar x,\bar t,\bar y,\bar s,\bar v,\bar w)$.
Because of \erf{a.3.4}, we find that $(\bar x,\bar t)\in B$ and $\bar v>\bar w$, and
we may assume by passing a subsequence of $\{\ga_j\}$ if necessary
that $(\hat x,\hat t),\,(\hat y,\hat s)\in B$, $v^{\ep,\eta}(\hat x,\hat t)>w(\hat y,\hat s)$.
By the viscosity properties of $v^{\ep,\eta}$ and $w$, we get for $\ga=\ga_j$,
\beq\label{a.3.5}\left\{\bald
&\ga^2(\hat t-\hat s)+\hH(\hat x, \ga(\hat x-\hat y),\hat t, v^{\ep,\eta}(\hat x,\hat t))+\gam e^{\gL \hat t}|\ga(\hat x-\hat y)|^2+(1-\eta)\gd\leq 0,
\\&\ga^2(\hat t-\hat s)+\hH(\hat y,\ga(\hat x-\hat y),\hat s, w(\hat y,\hat s))
+\gam e^{\gL\hat s}|\ga(\hat x-\hat y)|^2\geq 0.
\eald\right.
\eeq
Since $v^{\ep,\eta}(\hat x,\hat t)
>w(\hat y,\hat s)$, we have
\[
\hH(\hat x, \ga(\hat x-\hat y),\hat t, v^{\ep,\eta}(\hat x,\hat t))
\geq \hH(\hat x, \ga(\hat x-\hat y),\hat t, w(\hat x,\hat t)).
\]
Hence, from \erf{a.3.5} we obtain
\beq\label{a.3.6}\bald
&\tH(\hat x ,\ga(\hat x-\hat y),\hat t, w(\hat y,\hat s))
-\tH(\hat y,\ga(\hat x-\hat y),\hat s, w(\hat y,\hat s))
\\&\qquad +\gam(e^{\gL \hat t}-e^{\gL \hat s})|\ga(\hat x-\hat y|^2+(1-\eta)\gd \leq 0.
\eald
\eeq
Observe that
\[\bald
\tH(\hat x ,&\ga(\hat x-\hat y),\hat t, w(\hat y,\hat s))
\geq \tH(\hat y,\ga(\hat x-\hat y),\hat t, w(\hat y,\hat s))
-\go(|\hat x-\hat y|)
\\& \geq \tH(\hat y,\ga(\hat x-\hat y),\hat s, w(\hat y,\hat s))
-|e^{\gL\hat t}-e^{\gL\hat s}|
|H(\hat y,e^{\gL \hat t}\ga(\hat x-\hat y), e^{\gL\hat t}w(\hat y,\hat s))|
\\&\quad -\go(|e^{\gL\hat t}-e^{\gL\hat s}| (\ga|\hat x-\hat y|+|w(\hat y,\hat s)|))
-\go(|\hat x-\hat y|),
\eald
\]
where $\go$ denotes the modulus of continuity of the function $H$,
and moreover that
\begin{align*}
&|e^{\gL\hat t}-e^{\gL \hat s}|\leq \gL e^{\gL T}|\hat t-\hat s|\leq \gL
e^{\gL T}\ga|\hat t-\hat s|,
\\&\ga|\hat x-\hat y||\hat t-\hat s|\leq \ga(|\hat x-\hat y|^2+|\hat t-\hat s|^2)
\leq \ga |\hat x-\hat y|^2+\ga^2 |\hat t-\hat s|^2,
\\& |\hat t-\hat s|\ga^2|\hat x-\hat y|^2\leq \ga^2(|\hat x-\hat y|^4+|\hat t-\hat s|^2).
\end{align*}
Combine these observations with \erf{a.3.2}, we find that if $\ga=\ga_j$ is large enough, we have
\[\bald
&\tH(\hat x ,\ga(\hat x-\hat y),\hat t, w(\hat y,\hat s))
-\tH(\hat y,\ga(\hat x-\hat y),\hat s, w(\hat y,\hat s))
\\&\qquad +\gam(e^{\gL \hat t}-e^{\gL \hat s})|\ga(\hat x-\hat y|^2
>-(1-\eta)\gd.
\eald
\]
This contradicts \erf{a.3.6}, which completes the proof.
\eproof

\bproof[Proof of Theorem \ref{a.1}] We may assume that
$T<\infty$. Choose a constant $C_1>0$ so that
\[
|H(x,D\phi(x),\phi(x))|+|D\phi(x)|^2\leq C_1 \ \ \FORALL\ x\in M.
\]
Define the functions $f^\pm\in C^1(M\tim[0,T),\R)$ by
\[
f^+(x,t)=\phi(x)+C_1\gL^{-1}(e^{\gL t}-1) \ \ \AND f^-(x,t)=\phi(x)-C_1\gL^{-1}(e^{\gL t}-1).
\]
Choose a constant $C_2>0$ so that $-C_2\leq f^-\leq f^+\leq C_2$ on $M\tim[0,T)$.
We define the new Hamiltonians $\tH,\,\tH_k,\,\hH_k$ as follows. Set
\[\bald
\tH(x,p,u)&=H\left(x,p,\max\{-C_2,\min\{C_2,\,u\}\}\right) \ \ \FOR\ (x,p,u)\in\T^*M\tim\R,
\\ \tH_k(x,p,u)&=\min\{k,\tH(x,p,u)\} \ \ \ \ \FOR\ (x,p,u)\in T^*M\tim\R,\ k\in\N,
\\ \hH_k(x,p,u)&=\tH_k(x,p,u)+\tfrac 1k |p|^2 \ \ \ \ \FOR \ (x,p,u)\in T^*M\tim\R,\,k\in\N.
\eald
\]
It is easily seen that the functions $\tH(x,p,u),\, \tH_k(x,p,u),\,\hH_k(x,p,u)$ are
Lipschitz continuous in $u$, with $\gL$ as a Lipschitz bound, that for any
$f\in C^1(M\tim(0,T),\R)$, if $|f|\leq C_2$ on $T^*M\tim(0,T)$, then
$\tH(x,D_xf(x,t),f(x,t))=H(x,D_xf(x,t),f(x,t))$ for all $(x,t)\in M\tim(0,T)$, and that
$|\tH_k(x,D\phi(x),\phi(x))|\leq |\tH(x,D\phi(x),\phi(x))|$ \ and \
$|\hH_k(x,D\phi(x),\phi(x))|\leq C_1$ for all \ $x\in M$ and $k\in\N$.

Compute that
\[
f^+_t+H(x,D_xf^+,f^+)\geq C_1 e^{\gL t}+\hH_k(x,D_x\phi,\phi)-C_1(e^{\gL t}-1)\geq 0,
\]
to see that for any $k\in\N$, $f^+$ is a classical supersolution of
\beq\label{a.4}
u_t+\hH_k(x,D_xu,u)=0 \ \ \IN M\tim(0,T).
\eeq
Similarly, we find that for any $k\in\N$, $f^-$ is a classical subsolution of \erf{a.4}. Moreover, note that
$f^-(x,0)=f^+(x,0)=\phi(x)$ for all $x\in M$ and $f^-\leq f^+$ on $M\tim [0,T)$.
The Perron method yields a Crandall-Lions viscosity solution of \erf{a.4}. That is,
the formula
\[
u^k(x,t)=\sup\{f(x,t)\mid f\in\cS^-(\pl_t+\hH_k),\ f^-\leq f \leq f^+ \ \ \ON M\tim(0,T)\}
\]
gives a solution of \erf{a.4} such that $f^-\leq u^k\leq f^+$ on $M\tim(0,T)$.
Since $f^\pm(x,0)=\phi(x)$ for all $x\in M$, we may extend the domain of $u^k$
to $M\tim[0,T)$ so that
\[
u^k(x,0)=\phi(x)=
\lim_{M\tim(0,T)\ni (y,s)\to (x,0)}u(y,s) \ \ \FORALL \ x\in M.
\]
We note that, thanks to \erf{H2}, $\tH_k(x,p,u)=k$ if $|p|$ is {sufficiently} large
and, hence, the function $\tH_k$ is bounded and uniformly continuous on
$T^*M\tim\R$.  Since the upper and lower {semicontinuous} envelopes
$(u_k)^*$ and $(u_k)_*$ are respectively a viscosity subsolution and supersolution
of \erf{a.4}, we find by Lemma \ref{a.3} that $(u_k)^*\leq (u_k)_*$ on $M\tim[0,T)$,
which implies that $u_k\in C(M\tim[0,T),\R)$.

We now show that the family $\{u^k\}_{k\in\N}$ is equi-Lipschitz continuous
on $M\tim(0,T)$. For this, we show {first} that the functions $\hH_k$, with $k\in\N$,
are coercive uniformly in $k$. That is, for any $R>0$ there exists $Q>0$, chosen independently of $k$, such that for any $k\in\N$,
\beq\label{a.5}
\hH_k(x,p,u)>R \ \ \ \IF \ |p|>Q.
\eeq
Indeed, when $R>0$ is fixed, by \erf{H2} we may choose $Q\geq R$ so that
$\tH(x,p,u)>R$ if $|p|>Q$.  Using the inequality
\[
R<k+\tfrac 1 kR^2,
\]
we find that if $|p|>Q$, then
\[
\hH_k(x,p,u)\geq \min\{k+\tfrac 1k Q^2,\tH(x,p,u)\}>R.
\]
Hence, \erf{a.5} is valid.

Fix $h>0$ sufficiently small. Since $u^k\geq f^-$ on $M\tim[0,T)$, we find that
\[u^k(x,h)\geq \phi(x)-C_1\gL^{-1}(e^{\gL h}-1)\geq \phi(x)-C_1e^{\gL T}h \ \
\FORALL \ x\in M.
\]
Setting $v(x,t)=u^k(x,t)-C_1h e^{\gL(T+t)}$ for $(x,t)\in M\tim[0,T)$, we
easily observe that $v$ is a viscosity
subsolution of \erf{a.4}. We apply Lemma \ref{a.3},
to obtain the inequality $v(x,t)\leq u(x,h+t)$ for all $(x,t)\in M\tim[0,T-h)$, which shows that
\[
\liminf_{h\to 0+}\fr{u^k(x,t+h)-u^k(x,t)}{h}\geq -C_1e^{2\gL T}.
\]
This assures that for any $(x,t)\in M\tim(0,T)$ and $(p,q)\in D^+u^k(x,t)$,
$q\geq -C_1e^{2\gL T}$ and
\[
|q|+\hH_k(x,p,u^k(x,t)) \leq q+\hH_k(x,p,u^k(x,t))+2C_1e^{2\gL T}.
\]
Thus, $u^k$ is a viscosity subsolution of
\[
|u_t|+\hH_k(x,D_xu,u)-2C_1e^{2\gL T}=0 \ \ \IN M\tim(0,T),
\]
which, together with \erf{a.5}, ensures that for some constant $C_3>0$,
independent of $k$,
\[
|u^k_t|+|D_x u^k|\leq C_3 \ \ \IN\ M\tim(0,T) \ \text{  in the viscosity sense}.
\]
This shows (see \cite[Theorem I.14]{CL1983}, \cite[Proposition 1.14]{I2013})
that $\{u^k\}_{k\in\N}$ is equi-Lipschitz continuous on $M\tim[0,T)$.
Recalling that $f^-\leq u^k\leq f^+$ on $M\tim[0,T)$ for all $k\in\N$,
we find by the Ascoli-Arzela theorem that the family $\{u^k\}_{k\in\N}$ has
a subsequence, converging to some $u$ in $ C(M\tim[0,T),\R)$.
Since $\{\hH_k\}_{k\in\N}$ converges to $\tH$ in $C(T^*M\tim\R,\R)$,
we find that $u$ is a {viscosity} solution of
$u_t+\tH(x,D_xu,u)=0$ in $M\tim(0,T)$. It is obvious that $|u|\leq C_2$ on $M\tim[0,T)$, which implies that $u$ is a viscosity solution of $u_t+H(x,D_xu,u)=0$ in $M\tim(0,T)$, $u$ is Lipschitz continuous on $M\tim[0,T)$, and $u(x,0)=\phi(x)$ for all $x\in M$.
Thus, $u$ is a Lipschitz continuous solution of \erf{cp.1}.
The uniqueness assertion of the current theorem is a well-known result
and we do not repeat the standard proof here.  The uniqueness is also a consequence of Theorem \ref{cwch}. It follows as well from Lemma \ref{a.3},
once the Hamiltonian $H$ is replaced by a bounded and uniformly continuous function, which can be done based on the Lipschitz continuity of given solutions.
\eproof

\bye